\documentclass[journal]{IEEEtran}
\usepackage{cite}

\ifCLASSINFOpdf
  \usepackage[pdftex]{graphicx}
\else
\fi
\usepackage[cmex10]{amsmath}
\usepackage{algorithmic}
\usepackage[ruled, vlined]{algorithm2e}

\usepackage{array}

\usepackage{stfloats}
\usepackage{url}

\usepackage{color}

\usepackage[tight,footnotesize]{subfigure}
\usepackage{labtex}
\usepackage{enumerate}

\newcommand{\s}{\mathcal{S}}
\newtheorem{remark}{Remark}
\newcommand{\RHinf}{\mathcal{RH}_\infty}

\newcommand{\Set}[1]{\mathpzc{#1}}
\newcommand{\Sub}[3]{#1(#2, #3)}
\newcommand{\ttf}[1]{\boldsymbol{#1}}
\newcommand{\tf}[1]{\mathbf{#1}}

\newcommand{\C}{\mathcal C}
\newcommand{\K}{\mathcal{C}_\mathcal{K}}

\newcommand{\Ell}{\mathcal L}
\newcommand{\FT}{\mathcal F_T}
\newcommand{\Sother}{\mathcal{X}}

\newcommand{\coldec}[1]{#1^{(c)}}
\newcommand{\rowdec}[1]{#1^{(r)}}
\newcommand{\colps}{\Set c}
\newcommand{\rowps}{\Set r}
\newcommand{\rone}{\Set s}
\newcommand{\rtwo}{\Set t}

\begin{document}
%
\title{Separable and Localized System Level Synthesis for Large-Scale Systems}
%
%
%

\author{Yuh-Shyang~Wang,~\IEEEmembership{Member,~IEEE,}
        Nikolai~Matni,~\IEEEmembership{Member,~IEEE,}
        and~John~C.~Doyle
\thanks{The authors are with the department of Control and Dynamical Systems, California Institute of Technology, Pasadena, CA 91125, USA ({\tt\small \{yswang,nmatni,doyle\}@caltech.edu}).}%
\thanks{This research was in part supported by NSF NetSE, AFOSR, the Institute for Collaborative Biotechnologies through grant W911NF-09-0001 from the U.S. Army Research Office, and from MURIs ``Scalable, Data-Driven, and Provably-Correct Analysis of Networks'' (ONR) and ``Tools for the Analysis and Design of Complex Multi-Scale Networks'' (ARO). The content does not necessarily reflect the position or the policy of the Government, and no official endorsement should be inferred.}%
\thanks{Preliminary versions of this work have appeared in \cite{2014_Wang_CDC, 2015_Wang_LDKF, 2015_Wang_H2, 2015_Wang_Reg}.}
}

%
%

\markboth{Submitted to IEEE Transactions on Automatic Control,~Vol.~XX, No.~XX, January~2017}%
{Shell \MakeLowercase{\textit{et al.}}: Localized System Level Synthesis for Large-Scale Systems}
%



\maketitle

\begin{abstract}
A major challenge faced in the design of large-scale cyber-physical systems, such as power systems, the Internet of Things or intelligent transportation systems, is that traditional distributed optimal control methods do not scale gracefully, neither in controller synthesis nor in controller implementation, to systems composed of millions, billions or even trillions of interacting subsystems.  This paper shows that this challenge can now be addressed by leveraging the recently introduced System Level Approach (SLA) to controller synthesis.  In particular, in the context of the SLA, we define suitable notions of separability for control objective functions and system constraints such that the global optimization problem (or iterate update problems of a distributed optimization algorithm) can be decomposed into parallel subproblems.  We then further show that if additional locality (i.e., sparsity) constraints are imposed, then these subproblems can be solved using local models and local decision variables.  The SLA is essential to maintaining the convexity of the aforementioned problems under locality constraints.  As a consequence, the resulting synthesis methods have $O(1)$ complexity relative to the size of the global system.  We further show that many optimal control problems of interest, such as (localized) LQR and LQG, $\mathcal{H}_2$ optimal control with joint actuator and sensor regularization, and (localized) mixed $\mathcal{H}_2/\mathcal{L}_1$ optimal control problems, satisfy these notions of separability, and use these problems to explore tradeoffs in performance, actuator and sensing density, and average versus worst-case performance for a large-scale power inspired system.
\end{abstract}


%
\IEEEpeerreviewmaketitle

\section*{Preliminaries \& Notation}
We use lower and upper case Latin letters such as $x$ and $A$ to denote vectors and matrices, respectively, and lower and upper case boldface Latin letters such as $\tf x$ and $\tf G$ to denote signals and transfer matrices, respectively.  We use calligraphic letters such as $\s$ to denote sets.

We work with discrete time linear time invariant systems. We use standard definitions of the Hardy spaces $\mathcal{H}_2$ and $\mathcal{H}_\infty$, and denote their restriction to the set of real-rational proper transfer matrices by $\mathcal{RH}_2$ and $\RHinf$.  We use $G[i]$ to denote the $i$th spectral component of a transfer function $\tf G$, i.e., $\tf G(z) = \sum_{i=0}^{\infty} \frac{1}{z^i} G[i]$ for $| z | > 1$.  Finally, we use $\FT$ to denote the space of finite impulse response (FIR) transfer matrices with horizon $T$, i.e., $\FT := \{ \tf G \in \RHinf \, | \, \tf G = \sum_{i=0}^T\frac{1}{z^i}G[i]\}$.

Let $\mathbb{Z}^+$ be the set of all positive integers. We use calligraphic lower case letters such as $\Set r$ and $\Set c$ to denote subsets of $\mathbb{Z}^+$. 
Consider a transfer matrix $\ttf{\Phi}$ with $n_r$ rows and $n_c$ columns. Let $\rowps$ be a subset of $\{1, \dots, n_r \}$ and $\colps$ a subset of $\{1, \dots, n_c \}$. We use $\Sub{\ttf{\Phi}}{\rowps}{\colps}$ to denote the submatrix of $\ttf{\Phi}$ by selecting the rows according to the set $\rowps$ and columns according to the set $\colps$. 
We use the symbol $:$ to denote the set of all rows or all columns, i.e., we have $\ttf{\Phi} = \Sub{\ttf{\Phi}}{:}{:}$. 
Let $\{ \colps_1, \dots \colps_p \}$ be a partition of the set $\{1, \dots, n_c\}$. Then $\{ \Sub{\ttf{\Phi}}{:}{\colps_1}, \dots,  \Sub{\ttf{\Phi}}{:}{\colps_p} \}$ is a column-wise partition of the transfer matrix $\ttf{\Phi}$.
\section{Introduction}

\IEEEPARstart{L}{arge-scale} networked systems have emerged in extremely diverse application areas, with examples including the Internet of Things, the smart grid, automated highway systems, software-defined networks, and biological networks in science and medicine. The scale of these systems poses fundamental challenges to controller design: simple locally tuned controllers cannot guarantee optimal performance (or at times even stability), whereas a traditional centralized controller is neither scalable to compute nor physically implementable. Specifically, the synthesis of a centralized optimal controller requires solving a large-scale optimization problem that relies on the {global} plant model. In addition, information in the control network is assumed to be exchanged {instantaneously}. For large systems, the computational burden of computing a centralized controller, as well as the degradation in performance due to communication delays between sensors, actuators and controllers, make a centralized scheme unappealing.

The field of distributed (decentralized) optimal control developed to address these issues, with an initial focus placed on
explicitly incorporating communication constraints between sensors, actuators, and controllers into the control design process.
These communication delays are typically imposed as information sharing constraints on the set of admissible controllers in the resulting optimal control problem. 
The additional constraint makes the distributed optimal control problem significantly harder than its unconstrained centralized counterpart, and certain constrained optimal control problems are indeed NP-hard \cite{1968_Witsenhausen_counterexample, 1984_Tsitsiklis_NP_hard}. 
It has been shown that in the model-matching framework, distributed optimal control problems admit a convex reformulation in the Youla domain if \cite{2006_Rotkowitz_QI_TAC} and only if \cite{2014_Lessard_convexity} the information sharing constraint is quadratically invariant (QI) with respect to the plant.
With the identification of quadratic invariance as a means of convexifying the distributed optimal control problem, tractable solutions for various types of distributed objectives and constraints have been developed \cite{2012_Lessard_two_player, 2014_Lamperski_H2_journal, 2014_Lamperski_state, 2014_Sabau_QI, 2014_Matni_Hinf, 2014_Lessard_Hinf, 2013_Scherer_Hinf}.

As promising and impressive as all of the results have been, QI focuses primarily on identifying the tractability (e.g., the convexity) of the constrained optimal control problem, and places less emphasis on the scalability of synthesizing and implementing the resulting constrained optimal controller. 
In particular, for a strongly connected networked system,\footnote{A networked system is strongly connected if the state of any subsystem can eventually alter the state of all other subsystems.} quadratic invariance is only satisfied if each sub-controller communicates its locally acquired information with every other sub-controller in the system -- this density in communication has negative consequences on the scalability of both the synthesis and implementation of such distributed controllers.

To address these scalability issues, techniques based on regularization \cite{2011_Fardad_sparsity_promoting, 2013_Lin_sparsity_promoting}, convex approximation \cite{2014_Fazelnia_ODC, 2014_Sznaier, 2015_DJ}, and spatial truncation \cite{2009_Motee_spatial_truncation} have been used to find sparse (static) feedback controllers that are scalable to implement (i.e., local control actions can be computed using a local subset of global system state). 
These methods have been successful in extending the size of systems for which a distributed controller can be implemented and computed, but there is still a limit to their scalability as they rely on an underlying centralized synthesis procedure.

In our companion paper \cite{2015_PartI}, we proposed a \emph{System Level Approach} (SLA) to controller synthesis, and showed it to be a generalization of the distributed optimal control framework studied in the literature. 
The key idea is to directly design the entire response of the feedback system first, and then use this system response to construct an internally stabilizing controller that achieves the desired closed loop behavior. Specifically, the constrained optimal control problem is generalized to the following \emph{System Level Synthesis} (SLS) problem:
\begin{subequations} \label{eq:clock}
\begin{align}
\underset{\ttf{\Phi}}{\text{minimize  }} \quad & g(\ttf{\Phi}) \label{eq:clock1}\\
\text{subject to }  \quad & \ttf{\Phi} \text{ stable and achievable } \label{eq:clock2}\\
& \ttf{\Phi} \in \s,  \label{eq:clock3}
\end{align}
\end{subequations}
where $\ttf{\Phi}$ is the system response from disturbance to state and control action (formally defined in Section \ref{sec:prelim}), $g(\cdot)$ is a functional that quantifies the performance of the closed loop response, the first constraint ensures that an internally stabilzing controller exists that achieves the synthesized response $\ttf{\Phi}$, and $\s$ is a set. 
An important consequence of this approach to constrained controller synthesis is that for any system, convex constraints can always be used to impose sparse structure on the resulting controller.  If the resulting SLS problem is feasible, it follows that the resulting controller is scalable to implement, as this sparse structure in the controller implies that only local information (e.g., measurements and controller internal state) need to be collected and exchanged to compute local control laws.

In this paper, we show that the SLA developed in \cite{2015_PartI} further allows for constrained optimal controllers to be synthesized in a \emph{localized} way, i.e., by solving (iterate) sub-problems of size scaling as $O(1)$ relative to the dimension of the global system.  
Specifically, we first introduce the notion of \emph{separability} of the objective functional $g(\cdot)$ and the convex set constraint $\s$ such that the global optimization problem \eqref{eq:clock} can be \emph{solved} via parallel computation.  We then show that if additional locality (sparsity) constraints are imposed, then these subproblems can be solved in a localized way.  We show that such locality and separability conditions are satisfied by many problems of interest, such as the Localized LQR (LLQR) and Localized Distributed Kalman Filter (LDKF) \cite{2014_Wang_CDC, 2015_Wang_LDKF}.
We then introduce the weaker notion of partially separable objectives and constraints sets, and show that this allows for iterate subproblems of distributed optimization techniques such as the alternating direction method of multipliers (ADMM) \cite{BoydADMM} to be solved via parallel computation.  Similarly to the separable case, when additional locality constraints are imposed, the iterate subproblems can further be solved in a localized way.  We show that a large class of SLS problems are partially separable. Examples include the (localized) mixed $\mathcal{H}_2 / \mathcal{L}_1$ optimal control problem and the (localized) $\mathcal{H}_2$ optimal control problem with sensor and actuator regularization \cite{MC_CDC14, 2015_Wang_Reg}. 

The rest of the paper is organized as follows. In Section \ref{sec:prelim}, we formally define networked systems as a set of interacting discrete-time LTI systems, recall the distributed optimal control formulation and comment on its scalability limitations, and summarize the SLA to controller synthesis defined in our companion paper \cite{2015_PartI}.
In Section \ref{sec:sf}, we introduce the class of \emph{column/row-wise separable problems} as a special case of \eqref{eq:clock}, and show that such problems decompose into subproblems that can be solved in parallel.
We then extend the discussion to the general case in Section \ref{sec:of}, where we define the class of \emph{partially separable problems}.
When the problem is partially separable, we use distributed optimization combined with the techniques introduced in Section \ref{sec:sf} to solve \eqref{eq:clock} in a scalable way. We also give many examples of partially separable optimal control problems, and introduce some specializations and variants of the optimization algorithm in the same section.
Simulation results are shown in Section \ref{sec:simu} to demonstrate the effectiveness of our method. Specifically, we synthesize a localized $\mathcal{H}_2$ optimal controller for a power-inspired system with $12800$ states composed of dynamically coupled randomized heterogeneous subsystems, demonstrate the optimal co-design of localized $\mathcal{H}_2$ controller with sensor actuator placement, and solve a mixed $\mathcal{H}_2/\mathcal{L}_1$ optimal control problem. Conclusions and future research directions are summarized in Section \ref{sec:conclusion}.

\section{Preliminaries and Problem Statement} \label{sec:prelim}
We first introduce the networked system model that we consider in this paper. We then explain why the traditional distributed optimal control problem, as studied in the QI literature, does not scale gracefully to large systems. We then recall the SLA to controller synthesis defined and analyzed in \cite{2015_PartI}. This section ends with the formal problem statement considered in this paper.


\subsection{Interconnected System Model} \label{sec:challenge}
We consider a discrete-time LTI system of the form
\begin{eqnarray}
x[t+1] &=& A x[t] + B_2 u[t] + \delta_x [t] \nonumber\\
\bar{z}[t] &=& C_1x[t] + D_{12}u[t] \nonumber\\
y[t] &=& C_2 x[t] + \delta_y [t] \label{eq:dynamics} 
\end{eqnarray}
where $x$ is the global system state, $u$ is the global system input, $\bar{z}$ is the global controlled output, $y$ is the global measured output, and $\delta_x$ and $\delta_y$ are the global process and sensor disturbances, respectively.  We assume that these disturbances are generated by an underlying disturbance $w$ satisfying $\delta_x = B_1 w$ and $\delta_y = D_{21}w$, allowing us to concisely describe the system in transfer function form as 
\begin{equation}
\tf P = \left[ \begin{array}{c|cc} A & B_1 & B_2 \\ \hline C_1 & {0} & D_{12} \\ C_2 & D_{21} & {0} \end{array} \right] = \begin{bmatrix} \tf P_{11} & \tf P_{12} \\ \tf P_{21} & \tf P_{22} \end{bmatrix} \nonumber
\end{equation}
where $\tf P_{ij} = C_i(zI-A)^{-1}B_j + D_{ij}$.  We use $n_x, n_y$, and $n_u$ to denote the dimension of $x, y$, and $u$, respectively.

We further assume that the global system is composed of $n$ dynamically coupled subsystems that interact with each other according to an interaction graph $\mathcal{G}=(\mathcal{V},\mathcal{E})$. Here $\mathcal{V} = \{1, \dots, n\}$ denotes the set of subsystems. We denote by $x_i, u_i,$ and $y_i$ the state vector, control action, and measurement of subsystem $i$ -- we further assume that the corresponding state-space matrices \eqref{eq:dynamics} admit compatible block-wise partitions. The set $\mathcal{E} \subseteq \mathcal{V} \times \mathcal{V}$ encodes the physical interaction between these subsystems -- an edge $(i,j)$ is in $\mathcal{E}$ if and only if the state $x_j$ of subsystem $j$ directly affects the state $x_i$ of subsystem $i$.

We also assume that the controller, which is the mapping from measurements to control actions, is composed of physically distributed sub-controllers interconnected via a communication network.  This communication network defines the information exchange delays that are imposed on the set of admissible controllers, which for the purposes of this paper, are encoded via subspace constraints (cf. \cite{2012_Lessard_two_player, 2014_Lamperski_H2_journal, 2014_Lamperski_state, 2014_Sabau_QI, 2014_Matni_Hinf, 2014_Lessard_Hinf, 2013_Scherer_Hinf} for examples of how information exchange delays can be encoded via subspace constraints).

The objective of the optimal control problem (formally defined in the next subsection) is then to design a feedback strategy from measurement $\tf y$ to control action $\tf u$, subject to the aforementioned information exchange constraints, that minimizes the norm of the closed loop transfer matrix from the disturbance $\tf w$ to regulated output $\tf{\bar{z}}$.


\begin{example}[Chain topology]
\begin{figure}[h!]
     \centering
     \includegraphics[width=0.45\textwidth]{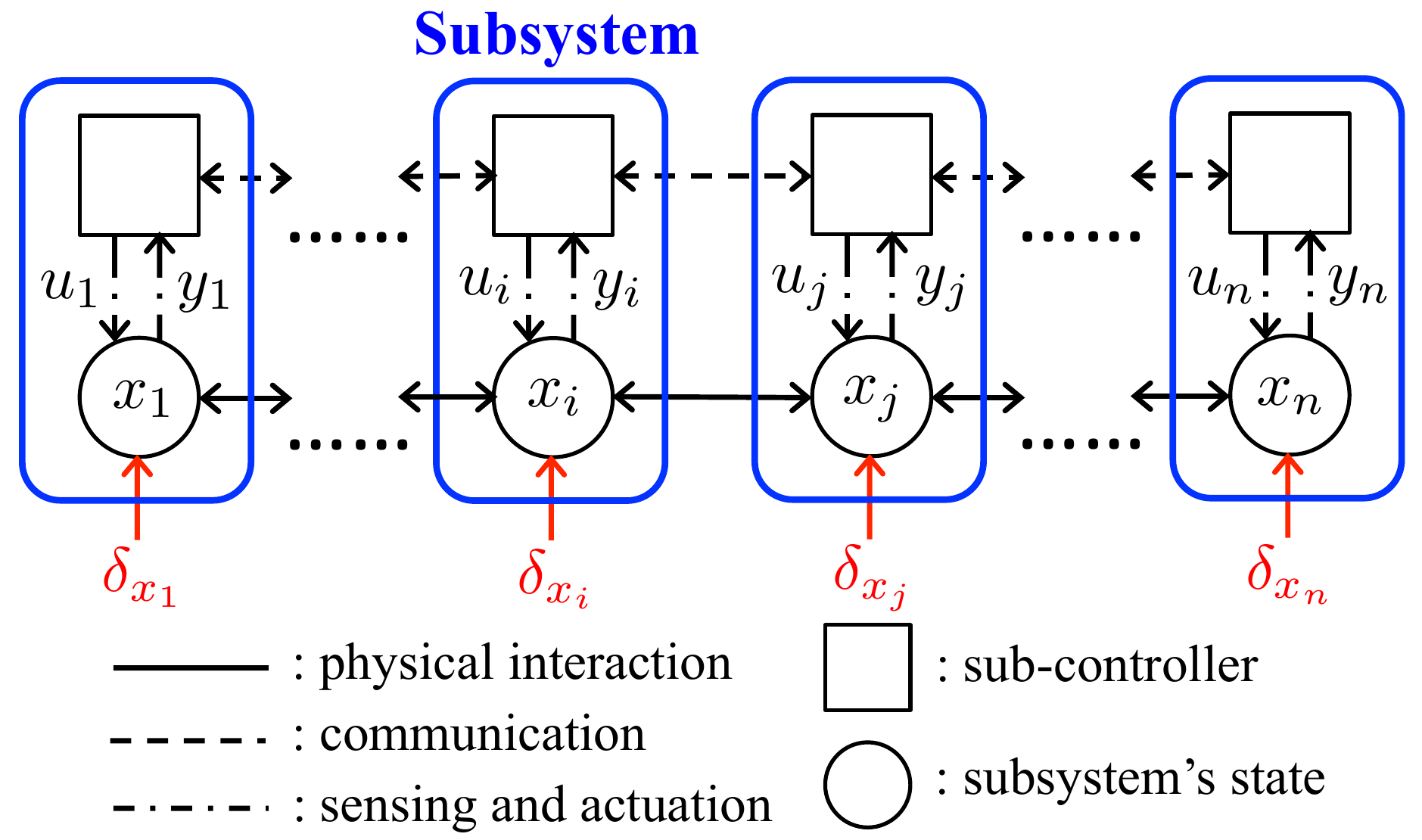}
     \caption{An Example of Interconnected System}
     \label{fig:interconnected}
\end{figure}
Figure \ref{fig:interconnected} shows an example of such an interconnected system.  The interaction graph of this system is a chain, and further, as illustrated, local control and disturbance inputs only directly affect local subsystems.  This allows us to write the dynamics of a subsystem $i$ as 
\begin{eqnarray}
x_i[t+1] &=& A_{ii} x_i[t] + \sum_{j \in \mathcal{N}_i} A_{ij}x_j[t] + B_{ii} u_i[t] + \delta_{x_i}[t] \nonumber\\
y_i[t] &=& C_{ii} x_i[t] + \delta_{y_i}[t], \label{eq:interconnected}
\end{eqnarray}
where $\mathcal{N}_i = \{j | (i,j) \in \mathcal{E}\}$ is the (incoming) neighbor set of subsystem $i$, $A_{ii}, A_{ij}, B_{ii}, C_{ii}$ are suitably defined sub-blocks of $A$, $B_2$ and $C_2$ matrices of compatible dimensions, and $\delta_{x_i}$ and $\delta_{y_i}$ are the process and sensor disturbances, respectively. 
\end{example}

\subsection{Distributed Optimal Control Framework}
Interconnected systems (such as that defined in equation \eqref{eq:interconnected}) pose fundamental challenges to controller synthesis: imposing information exchange constraints between sub-controllers can make the distributed optimal control problem non-convex, and even when these problems are tractable, the synthesis and implementation of the resulting optimal controller often does not scale well to large systems.  Here we briefly recall the distributed optimal control framework, and illustrate why scalability issues may arise when dealing with large systems.

Consider a dynamic output feedback control law $\tf u = \tf K \tf y$. The distributed optimal control problem of minimizing the norm of the closed loop transfer matrix from external disturbance $\tf w$ to regulated output $ \tf{\bar{z}}$ \cite{2006_Rotkowitz_QI_TAC, 2012_Lessard_two_player,2014_Lamperski_H2_journal, 2013_Scherer_Hinf, 2014_Lessard_Hinf, 2014_Matni_Hinf} is commonly formulated as
\begin{subequations}\label{eq:trad}
\begin{align}
\underset{\tf K}{\text{minimize  }} \quad & ||\tf P_{11} + \tf P_{12} \tf K(I-\tf P_{22}\tf K)^{-1} \tf P_{21}\| \label{eq:trad_obj}\\
\text{subject to  } \quad & \tf K \text{ internally stabilizes } \tf P \label{eq:trad_stable}\\ 
& \tf K \in \K \label{eq:trad_info}
\end{align}
\end{subequations}
where the constraint $\tf K \in \K$ enforces information sharing constraints between the sub-controllers. 
When $\K$ is a subspace constraint, it has been shown that \eqref{eq:trad} is convex if \cite{2006_Rotkowitz_QI_TAC} and only if \cite{2014_Lessard_convexity} $\K$ is quadratically invariant (QI) with respect to $\tf P_{22}$, i.e., $\tf K \tf P_{22} \tf K \in \K$ for all $\tf K \in \K$. 
Informally, quadratic invariance holds when sub-controllers are able to share information with each other at least as quickly as their control actions propagate through the plant \cite{2010_Rotkowitz_QI_delay}. 


When the interaction graph of the system \eqref{eq:dynamics} is strongly connected, the transfer matrix $\tf P_{22} = C_2 (zI-A)^{-1} B_2$ is generically dense, regardless as to whether the matrices $(A,B_2,C_2)$ are structured or not. In such cases, the conditions in \cite{2006_Rotkowitz_QI_TAC} imply that any sparsity constraint $\K$ imposed on $\tf K$ in \eqref{eq:trad} does not satisfy the QI condition, and therefore \eqref{eq:trad} cannot be reformulated in a convex manner in the Youla domain \cite{2014_Lessard_convexity}.
In other words, for a strongly connected system, the distributed optimal control problem is convex only when the transfer matrix $\tf K$ is dense, i.e., only when each local measurement $\tf y_i$ is shared among {all} sub-controllers $\tf u_j$ in the network.   For such systems, neither controller synthesis nor implementation scales gracefully to large systems.
Although recent methods based on convex relaxations \cite{2014_Fazelnia_ODC} can be used to solve certain cases of the non-convex optimal control problem \eqref{eq:trad} approximately for a sparse constraint set $\K$, the underlying synthesis problem is still large-scale and does not admit a scalable reformulation. The need to address scalability, both in the {synthesis} and {implementation} of a controller, is the driving motivation behind the SLA framework, which we recall in the next subsection.

\subsection{The System Level Approach}
In our companion paper \cite{2015_PartI}, we define and analyze the SLA to controller synthesis, defined in terms of three ``system level'' (SL) components: SL Parameterizations (SLPs), SL Constraints (SLCs) and SL Synthesis (SLS) problems.  We showed that the SLS problem is a substantial generalization of the distributed optimal control problem \eqref{eq:trad}.  We summarize the key results that we build upon here: conceptually, the key idea is to reformulate the optimal control problem in terms of the closed loop \emph{system response}, as opposed to the controller itself.

For an LTI system with dynamics given by \eqref{eq:dynamics}, we define a system response $\{\tf R, \tf M, \tf N, \tf L\}$ to be the maps satisfying
\begin{equation}
\begin{bmatrix} \tf x \\ \tf u \end{bmatrix} = \begin{bmatrix} \tf R & \tf N \\ \tf M & \tf L \end{bmatrix} \begin{bmatrix} \ttf{\delta_x} \\ \ttf{\delta_y} \end{bmatrix}. \label{eq:desired}
\end{equation}
We call a system response $\{\tf R, \tf M, \tf N, \tf L\}$ \emph{stable and achievable} with respect to a plant $\tf P$ if there exists an internally stabilizing controller $\tf K$ such that the control rule $\tf u = \tf K \tf y$ leads to closed loop behavior consistent with \eqref{eq:desired}. We showed in \cite{2015_PartI} that the parameterization of all stable achievable system responses $\{\tf R, \tf M, \tf N, \tf L\}$ is defined by the following affine space:
\begin{subequations} \label{eq:achievable}
\begin{align}
\begin{bmatrix} zI - A & -B_2 \end{bmatrix}
\begin{bmatrix} \tf R & \tf N \\ \tf M & \tf L \end{bmatrix} &= 
\begin{bmatrix} I & 0 \end{bmatrix} \label{eq:main_1} \\
\begin{bmatrix} \tf R & \tf N \\ \tf M & \tf L \end{bmatrix}
\begin{bmatrix} zI - A \\ -C_2 \end{bmatrix} &= 
\begin{bmatrix} I \\ 0 \end{bmatrix} \label{eq:main_2} \\
\tf R, \tf M, \tf N \in \frac{1}{z} \mathcal{RH}_{\infty}, & \quad \tf L \in \mathcal{RH}_{\infty}. \label{eq:main_stable}
\end{align}
\end{subequations}
The parameterization of all internally stabilizing controllers is further given by the following theorem.
\begin{theorem}[Theorem 2 in \cite{2015_PartI}] \label{thm:1} 
Suppose that a set of transfer matrices $(\tf R, \tf M, \tf N, \tf L)$ satisfy the constraints \eqref{eq:main_1} - \eqref{eq:main_stable}. Then an internally stabilizing controller yielding the desired system response \eqref{eq:desired} can be implemented as
\begin{subequations}\label{eq:implementation}
\begin{align}
z \ttf{\beta} &= \tilde{\tf R}^+ \ttf{\beta} + \tilde{\tf N} \tf y \label{eq:imple1}\\
\tf u &=  \tilde{\tf M} \ttf{\beta} + \tf L \tf y \label{eq:imple2}
\end{align}
\end{subequations}
where $\tilde{\tf R}^+ = z(I - z\tf R)$, $\tilde{\tf N} = -z\tf N$, $\tilde{\tf M} = z\tf M$, $\tf L$ are in $\RHinf$. Further, the solutions of \eqref{eq:main_1} - \eqref{eq:main_stable} with the implementation \eqref{eq:implementation} parameterize all internally stabilizing controllers. 
\end{theorem}
\begin{remark}
Equation \eqref{eq:implementation} can be considered as a natural extension of the state space realization of a controller $\tf K$, which is given by
\begin{eqnarray}
z \ttf{\xi} &=& A_K \ttf{\xi} + B_K \tf y \nonumber\\
\tf u &=& C_K \ttf{\xi} + D_K \tf y. \label{eq:sss}
\end{eqnarray}
Here we allow the state-space matrices $A_K, B_K, C_K, D_K$ of the controller, as specified in \eqref{eq:sss}, to instead be stable proper transfer matrices $\tf{\tilde{R}^{+}}, \tf{\tilde{M}}, \tf{\tilde{N}}, \tf L$ in \eqref{eq:implementation}. As a result, we call the variable $\ttf{\beta}$ in \eqref{eq:imple2} the \emph{controller state}.
\end{remark}

With the parameterization of all stable achievable system response in \eqref{eq:main_1} - \eqref{eq:main_stable}, we now introduce a convex objective functional $g(\cdot)$ and an additional convex set constraint $\s$ imposed on the system response. We call $g(\cdot)$ a system level objective (SLO) and $\s$ a system level constraint (SLC). The complete form of the SLS problem \eqref{eq:clock} is then given by
\begin{subequations} \label{eq:LOC}
\begin{align}
\underset{\{\tf R, \tf M, \tf N, \tf L\}}{\text{minimize  }} \quad & g(\tf R,\tf M,\tf N,\tf L) \label{eq:main_obj}\\
\text{subject to }  \quad & \eqref{eq:main_1} - \eqref{eq:main_stable} \label{eq:main_affine}\\
& \begin{bmatrix} \tf R & \tf N \\ \tf M & \tf L \end{bmatrix} \in \s. \label{eq:main_S}
\end{align}
\end{subequations}
We showed in \cite{2015_PartI} that the SLS problem \eqref{eq:LOC} is a generalization of the distributed optimal control problem \eqref{eq:trad} in the following sense: (i) any quadratically invariant subspace constraint imposed on the control problem \eqref{eq:trad_info} can be formulated as a special case of a convex SLC in \eqref{eq:main_S}, and (ii) the objective \eqref{eq:trad_obj}, which can be written as
\begin{equation}
||\begin{bmatrix} C_1 & D_{12} \end{bmatrix} \begin{bmatrix} \tf R & \tf N \\ \tf M & \tf L \end{bmatrix} \begin{bmatrix} B_1 \\ D_{21} \end{bmatrix} ||, \label{eq:trad_obj2}
\end{equation}
is a special case of a convex SLO in \eqref{eq:main_obj}. In particular, the SLS problem defines the broadest known class of constrained optimal control problems that can be solved using convex programming \cite{2015_PartI}.

\subsection{Localized Implementation} \label{sec:imp}
In our prior work \cite{2014_Wang_CDC, 2015_Wang_H2, 2015_PartI}, we showed how to use the SLC \eqref{eq:main_S} to design a controller that is scalable to implement. 
We first express the SLC set $\s$ as the intersection of three convex set components: a locality (sparsity) constraint $\Ell$, a finite impulse response (FIR) constraint $\FT$, and an arbitrary convex set component $\Sother$, i.e., $\s = \Ell \cap \FT \cap \Sother$. The locality constraint $\Ell$ imposed on a transfer matrix $\tf G$ is a collection of sparsity constraints of the form $\tf G_{ij} = 0$ for some $i$ and $j$. The constraint $\FT$ restricts the optimization variables to be a finite impulse response transfer matrix of horizon $T$, thus making optimization problem \eqref{eq:LOC} finite dimensional. The constraint $\Sother$ includes any other convex constraint imposed on the system, such as communication constraints, bounds on system norms, etc. This leads to a \emph{localized SLS problem} given by
\begin{subequations} \label{eq:localized}
\begin{align}
\underset{\{\tf R, \tf M, \tf N, \tf L\}}{\text{minimize  }} \quad & g(\tf R,\tf M,\tf N,\tf L) \label{eq:localized_slo}\\
\text{subject to }  \quad & \eqref{eq:main_1} - \eqref{eq:main_stable} \label{eq:localized_slp}\\
& \begin{bmatrix} \tf R & \tf N \\ \tf M & \tf L \end{bmatrix} \in \Ell \cap \FT \cap \Sother. \label{eq:localized_slc}
\end{align}
\end{subequations}

To design a controller that is scalable to implement, i.e., one in which each sub-controller needs to collect information from $O(1)$ other sub-controllers to compute its local control action, is to impose sparsity constraints on the system response $(\tf R, \tf M, \tf N, \tf L)$ through the locality constraint $\Ell$. 
Theorem \ref{thm:1} makes clear that the structure of the system response $(\tf R, \tf M, \tf N, \tf L)$ translates directly to the structure, and hence the implementation complexity, of the resulting controller \eqref{eq:implementation}. 
For instance, if each row of these transfer matrices has a small number of nonzero elements, then each sub-controller only needs to collect a small number of measurements and controller states to compute its control action \cite{2014_Wang_CDC, 2015_Wang_H2, 2015_PartI}. 
This holds true for arbitrary convex objective functions $g(\cdot)$ and arbitrary convex constraints $\Sother$. 

\begin{remark}
For a state feedback problem ($C_2 = I$ and $D_{21} = 0$), we showed in \cite{2015_PartI} that the localized SLS problem \eqref{eq:localized} can be simplified to
\begin{eqnarray}
\underset{\{\tf R, \tf M\}}{\text{minimize    }} && g(\tf R, \tf M) \nonumber\\
\text{subject to } && \begin{bmatrix} zI - A & -B_2 \end{bmatrix} \begin{bmatrix} \tf R \\ \tf M \end{bmatrix} = I \nonumber \\
&& \begin{bmatrix} \tf R \\ \tf M \end{bmatrix} \in \Ell \cap \FT \cap \Sother. \label{eq:sf_lop}
\end{eqnarray}
In addition, the controller achieving the desired system response can be implemented directly using the transfer matrices $\tf R$ and $\tf M$ \cite{2015_PartI}.
\end{remark}
\begin{remark}
It should be noted that although the optimization problems \eqref{eq:localized} and \eqref{eq:sf_lop} are convex for arbitrary locality constraints $\Ell$, they are not necessarily feasible.  Hence in the SLA framework, the burden is shifted from verifying the convexity of a structured optimal control problem to verifying its feasibility.
In \cite{2015_Wang_Reg}, we showed that a necessary condition for the existence of a localized (sparse) system response is that the communication speed between sub-controllers is faster than the speed of disturbance propagation in the plant.
This condition can be more restrictive than the delay conditions \cite{2010_Rotkowitz_QI_delay} that must be satisfied for a system to be QI . However, when such a localized system response exists, we only require fast but \emph{local} communication to implement the controller.
\end{remark}

\subsection{Problem Statement}
In previous work \cite{2015_PartI,2014_Wang_CDC, 2015_Wang_LDKF, 2015_Wang_H2, 2015_Wang_Reg}, we showed that by imposing additional constraints $\Ell$ and $\FT$ in \eqref{eq:localized}, we can design a localized optimal controller that is scalable to implement. In this paper, we further show that the localized SLS problem \eqref{eq:localized} can be {solved} in a localized and scalable way if the SLO \eqref{eq:localized_slo} and the SLC \eqref{eq:localized_slc} satisfy certain \emph{separability} properties. Specifically, when the localized SLS problem \eqref{eq:localized} is separable and the locality constraint $\Ell$ is suitably specified, the global problem \eqref{eq:localized} can be decomposed into parallel local subproblems. The complexity of synthesizing each sub-controller in the network is then independent of the size of the global system. 

\section{Column-wise Separable SLS Problems} \label{sec:sf}

In this section, we consider the state feedback localized SLS problem \eqref{eq:sf_lop}, which is a special case of \eqref{eq:localized}.
We begin by defining the notion of a column-wise separable optimization problem. We then show that under suitable assumptions on the objective function and the additional constraint set $\mathcal{S}$, the state-feedback SLS problem \eqref{eq:sf_lop} satisfies these conditions. We then give examples of optimal control problems, including the (localized) LQR \cite{2014_Wang_CDC} problem, that belong to the class of column-wise separable problems.
In Section \ref{sec:reduction-p}, we propose a dimension reduction algorithm to further reduce the complexity of each column-wise subproblem from global to local scale. 
This section ends with an overview of row-wise separable SLS problems that are naturally seen as dual to their column-wise separable counterparts.

\subsection{Column-wise Separable Problems} \label{sec:gopt}
We begin with a generic optimization problem given by
\begin{subequations} \label{eq:gopt}
\begin{align} 
\underset{\ttf{\Phi}}{\text{minimize }} \quad & g(\ttf{\Phi}) \label{eq:gopt-1}\\
\text{subject to } \quad & \ttf{\Phi} \in \s, \label{eq:gopt-2}
\end{align}
\end{subequations}
where $\ttf{\Phi}$ is a $m \times n$ transfer matrix, $g(\cdot)$ a functional objective, and $\s$ a set constraint. Our goal is to exploit the structure of \eqref{eq:gopt} to solve the optimization problem in a computationally efficient way. Recall that we use calligraphic lower case letters such as $\Set c$ to denote subsets of positive integer. Let $\{ \colps_1, \dots, \colps_p\}$ be a partition of the set $\{1, \dots, n\}$. The optimization variable $\ttf{\Phi}$ can then be partitioned column-wise as $\{ \Sub{\ttf{\Phi}}{ : }{\colps_1}, \dots, \Sub{\ttf{\Phi}}{ : }{\colps_p} \}$. With the column-wise partition of the optimization variable, we can then define the column-wise separability of the functional objective \eqref{eq:gopt-1} and the set constraint \eqref{eq:gopt-2} as follows.

\begin{definition} \label{dfn:dd}
The functional objective $g(\ttf{\Phi})$ in \eqref{eq:gopt-1} is said to be \emph{column-wise separable} with respect to the column-wise partition $\{ \colps_1, \dots, \colps_p\}$ if
\begin{equation}
g(\ttf{\Phi}) = \sum_{j=1}^p g_j(\Sub{\ttf{\Phi}}{ : }{\colps_j}) \label{eq:gopt-obj}
\end{equation}
for some functionals $g_j(\cdot)$ for $j = 1, \dots, p$. 
\end{definition}

We now give a few examples of column-wise separable objectives.
\begin{example} [Frobenius norm]
The square of the Frobenius norm of a $m \times n$ matrix $\Phi$ is given by
\begin{equation}
\| \Phi \|_{\mathcal{F}}^2 = \sum_{i=1}^m \sum_{j=1}^n \Phi_{ij}^2. \nonumber
\end{equation}
This objective function is column-wise separable with respect to arbitrary column-wise partition.
\end{example}
\begin{example} [$\mathcal{H}_2$ norm]
The square of the $\mathcal{H}_2$ norm of a transfer matrix $\ttf{\Phi}$ is given by
\begin{equation}
\| \ttf{\Phi} \|_{\mathcal{H}_2}^2 = \sum_{t=0}^\infty \| \Phi[t] \|_{\mathcal{F}}^2, \label{eq:h2_comp}
\end{equation}
which is column-wise separable with respect to arbitrary column-wise partition.
\end{example}

\begin{definition} \label{dfn:s}
The set constraint $\s$ in \eqref{eq:gopt-2} is said to be \emph{column-wise separable} with respect to the column-wise partition $\{ \colps_1, \dots, \colps_p\}$ if the following condition is satisfied:
\begin{equation}
\ttf{\Phi} \in \s \quad \text{if and only if} \quad \Sub{\ttf{\Phi}}{ : }{\colps_j} \in \s_j \text{  for  } j = 1, \dots, p 
\label{eq:gopt-decom_scon}
\end{equation}
for some sets $\s_j$ for $j = 1, \dots, p$.
\end{definition}

%
\begin{example} [Affine Subspace]
The affine subspace constraint
\begin{equation}
\tf G \ttf{\Phi} = \tf H \nonumber
\end{equation}
is column-wise separable with respect to arbitrary column-wise partition. Specifically, we have
\begin{equation}
\tf G \Sub{\ttf{\Phi}}{ : }{\colps_j} = \Sub{\tf{H}}{ : }{\colps_j} \nonumber
\end{equation}
for $\colps_j$ any subset of $\{1, \dots, n\}$.
\end{example}
\begin{example} [Locality and FIR Constraints]
The locality and FIR constraints introduced in Section \ref{sec:imp}
\begin{equation}
\ttf{\Phi} \in \Ell \cap \FT \nonumber
\end{equation}
are column-wise separable with respect to arbitrary column-wise partition.  This follows from the fact that both locality and FIR constraints can be encoded via sparsity structure: the resulting linear subspace constraint is trivially column-wise separable.
\end{example}

Assume now that the objective function \eqref{eq:gopt-1} and the set constraint \eqref{eq:gopt-2} are both column-wise separable with respect to a column-wise partition $\{ \colps_1, \dots, \colps_p\}$. In this case, we say that optimization problem \eqref{eq:gopt} is a \emph{column-wise separable problem}. Specifically, \eqref{eq:gopt} can be partitioned into $p$ parallel subproblems as
\begin{subequations} \label{eq:gopt-decom}
\begin{align} 
\underset{\Sub{\ttf{\Phi}}{ : }{\colps_j}}{\text{minimize }} \quad & g_j(\Sub{\ttf{\Phi}}{ : }{\colps_j}) \label{eq:gopt-decom-1} \\
\text{subject to } \quad & \Sub{\ttf{\Phi}}{ : }{\colps_j} \in \s_j \label{eq:gopt-decom-3}
\end{align}
\end{subequations}
for $j = 1, \dots, p$. 

\subsection{Column-wise Separable SLS Problems} \label{sec:CWSP}

We now specialize our discussion to the localized SLS problem \eqref{eq:sf_lop}. 
We use
\begin{equation}
\ttf{\Phi} = \begin{bmatrix} \tf R \\ \tf M \end{bmatrix} \nonumber
\end{equation}
to represent the system response we want to optimize for, and we denote $\tf Z_{AB}$ the transfer matrix $\begin{bmatrix} zI-A & -B_2 \end{bmatrix}$. The state feedback localized SLS problem \eqref{eq:sf_lop} can then be written as 
\begin{subequations} \label{eq:sf_lop2}
\begin{align} 
\underset{\ttf{\Phi}}{\text{minimize }} \quad & g(\ttf{\Phi}) \label{eq:lop-1}\\
\text{subject to } \quad & \tf Z_{AB} \ttf{\Phi} = I \label{eq:lop-2} \\
& \ttf{\Phi} \in \s, \label{eq:lop-3}
\end{align}
\end{subequations}
with $\s = \Ell \cap \FT \cap \Sother$. Note that the affine subspace constraint \eqref{eq:lop-2} is column-wise separable with respect to any column-wise partition, and thus the overall problem is column-wise separable with respect to a given column-wise partition if the SLO \eqref{eq:lop-1} and SLC \eqref{eq:lop-3} are also column-wise separable with respect to that partition.
\begin{remark}
Recall that the locality constraint $\Ell$ and the FIR constraint $\FT$ are column-wise separable with respect to any column-wise partition. Therefore, the column-wise separability of the SLC $\s = \Ell \cap \FT \cap \Sother$ in \eqref{eq:sf_lop2} is determined by the column-wise separability of the constraint $\Sother$. If $\s$ is column-wise separable, then we can express the set constraint $\s_j$ in \eqref{eq:gopt-decom_scon} as $\s_j = \Sub{\Ell}{:}{\colps_j} \cap \FT \cap \Sother_j$ for some $\Sother_j$ for each $j$.
\end{remark}

Assume that the SLO \eqref{eq:lop-1} and the SLC \eqref{eq:lop-3} are both column-wise separable with respect to a column-wise partition $\{ \colps_1, \dots, \colps_p\}$. In this case, we say that the state feedback localized SLS problem \eqref{eq:sf_lop2} is a \emph{column-wise separable SLS problem}. Specifically, \eqref{eq:sf_lop2} can be partitioned into $p$ parallel subproblems as
\begin{subequations} \label{eq:decom}
\begin{align} 
\underset{\Sub{\ttf{\Phi}}{ : }{\colps_j}}{\text{minimize }} \quad & g_j(\Sub{\ttf{\Phi}}{ : }{\colps_j}) \label{eq:decom-1}\\
\text{subject to } \quad & \tf Z_{AB} \Sub{\ttf{\Phi}}{ : }{\colps_j} = \Sub{I}{ : }{\colps_j} \label{eq:decom-2} \\
& \Sub{\ttf{\Phi}}{ : }{\colps_j} \in \Sub{\Ell}{:}{\colps_j} \cap \FT \cap \Sother_j \label{eq:decom-3}
\end{align}
\end{subequations}
for $j = 1, \dots, p$. 

Here we give a few examples of column-wise separable SLS problems. 
\begin{example}[LLQR with uncorrelated disturbance]
In \cite{2014_Wang_CDC}, we formulate the LLQR problem with uncorrelated process noise (i.e., $B_1 = I$) as
\begin{subequations} \label{eq:llqr}
\begin{align}
\underset{\{\tf R, \tf M\}}{\text{minimize}} \quad & || \begin{bmatrix} C_1 & D_{12} \end{bmatrix} \begin{bmatrix} \tf R \\ \tf M \end{bmatrix} ||_{\mathcal{H}_2}^2 \label{eq:llqr-1}\\
\text{subject to} \quad & \begin{bmatrix} zI - A & -B_2 \end{bmatrix} \begin{bmatrix} \tf R \\ \tf M \end{bmatrix} = I \label{eq:llqr-2} \\
& \begin{bmatrix} \tf R \\ \tf M \end{bmatrix} \in \Ell \cap \FT \cap \frac{1}{z} \mathcal{RH}_\infty. \label{eq:llqr-3}
\end{align}
\end{subequations}
The SLO \eqref{eq:llqr-1} and the SLC \eqref{eq:llqr-3} are both column-wise separable with respect to arbitrary column-wise partition. The separable property of the SLO is implied by the separability of the $\mathcal{H}_2$ norm and the assumption that the process noise is pairwise uncorrelated. The separability of the constraints \eqref{eq:llqr-2} - \eqref{eq:llqr-3} is follows from the above discussion pertaining to affine subspaces. The physical interpretation of the column-wise separable property is that we can analyze the system response of each local process disturbance $\ttf{\delta_{x_j}}$ in an independent and parallel way, and then exploit the superposition principle satisfied by LTI systems to reconstruct the full solution to the LLQR problem.  We also note that removing the locality and FIR constraints does not affect column-wise separability, and hence the standard LQR problem with uncorrelated noise is also column-wise separable.
\end{example}
\begin{example}[LLQR with locally correlated disturbances] \label{ex:perm}
Building on the previous example, we now assume that the covariance of the process noise $\delta_x$ is given by $B_1^\top B_1$ for some matrix $B_1$, i.e., $\delta_x \sim \mathcal{N}(0,B_1^\top B_1)$. 
In addition, suppose that there exists a permutation matrix $\Pi$ such that the matrix $\Pi B_1$ is block diagonal. This happens when the global noise vector $\ttf{\delta_x}$ can be partitioned into uncorrelated subsets. The SLO for this LLQR problem can then be written as
\begin{equation}
|| \begin{bmatrix} C_1 & D_{12} \end{bmatrix} \ttf{\Phi} B_1 ||_{\mathcal{H}_2}^2  =  || \begin{bmatrix} C_1 & D_{12} \end{bmatrix} \ttf{\Phi} \Pi^\top \Pi B_1 ||_{\mathcal{H}_2}^2. \nonumber
\end{equation}
Note that $\ttf{\Phi} \Pi^\top$ is a column-wise permutation of the optimization variable. We can define a column-wise partition on $\ttf{\Phi} \Pi^\top$ according to the block diagonal structure of the matrix $\Pi B_1$ to decompose the SLO in a column-wise manner. The physical meaning is that we can analyze the system response of each uncorrelated subset of the noise vector in an independent and parallel way.
\end{example}
\begin{example}[Element-wise $\ell_1$] \label{ex:l1}
In the LLQR examples above, we rely on the separable property of the $\mathcal{H}_2$ norm, as shown in \eqref{eq:h2_comp}. 
Motivated by the $\mathcal{H}_2$ norm, we define the element-wise $\ell_1$ norm (denoted by $e_1$) of a transfer matrix $\tf G \in \FT$ as
\begin{equation}
|| \tf G ||_{e_1}  =  \sum_i \sum_j \sum_{t=0}^T | G_{ij}[t] |. \nonumber
\end{equation}
The previous examples still hold if we change the square of the $\mathcal{H}_2$ norm to the element-wise $\ell_1$ norm.
\end{example}
\subsection{Dimension Reduction Algorithm} \label{sec:reduction-p}
In this subsection, we discuss how the locality constraint $\Sub{\Ell}{:}{\colps_j}$ in \eqref{eq:decom-3} can be exploited to reduce the dimension of each subproblem \eqref{eq:decom} from global to local scale, thus yielding a localized synthesis procedure.
The general and detailed dimension reduction algorithm for arbitrary locality constraints can be found in Appendix \ref{sec:reduction}. Here we rather highlight the main results and consequences of this approach through the use of some examples.

In particular, in Appendix \ref{sec:reduction} we show that the optimization subproblem \eqref{eq:decom} is equivalent to
\begin{subequations} \label{eq:decom40}
\begin{align} 
\underset{\Sub{\ttf{\Phi}}{\rone_j}{\colps_j}}{\text{minimize }} \quad & \bar{g}_j(\Sub{\ttf{\Phi}}{\rone_j}{\colps_j}) \label{eq:decom40-1}\\
\text{subject to } \quad & \Sub{\tf Z_{AB}}{\rtwo_j}{\rone_j} \Sub{\ttf{\Phi}}{\rone_j}{\colps_j} = \Sub{I}{\rtwo_j}{\colps_j} \label{eq:decom40-2} \\
& \Sub{\ttf{\Phi}}{\rone_j}{\colps_j} \in \Sub{\Ell}{\rone_j}{\colps_j} \cap \FT \cap \bar{\Sother}_j \label{eq:decom40-3}
\end{align}
\end{subequations}
where $\rone_j$ and $\rtwo_j$ are sets of positive integers, and $\bar{g}_j$ and $\bar{\Sother}_j$ are objective function and constraint set restricted to the reduced dimensional space, respectively. 
Roughly speaking, the set $\rone_j$ is the collection of optimization variables contained within the localized region specified by $\Sub{\Ell}{:}{\colps_j}$, and the set $\rtwo_j$ is the collection of states that are directly affected by the optimization variables in $\rone_j$.
The complexity of solving \eqref{eq:decom40} is determined by the cardinality of the sets $\colps_j$, $\rone_j$, and $\rtwo_j$, which in turn are determined by the number of nonzero elements allowed by the locality constraint and the structure of the system state-space matrices \eqref{eq:dynamics}. For instance, the cardinality of the set $\rone_j$ is equal to the number of nonzero rows of the locality constraint $\Sub{\Ell}{:}{\colps_j}$.
When the locality constraint and the system matrices are suitably sparse, it is possible to make the size of these sets much smaller than the size of the global network. In this case, the global optimization subproblem \eqref{eq:decom} reduces to a local optimization subproblem \eqref{eq:decom40} which depends on the local plant model $\Sub{\tf Z_{AB}}{\rtwo_j}{\rone_j}$ only. 

\begin{figure}[h!]
     \centering
     \includegraphics[width=0.35\textwidth]{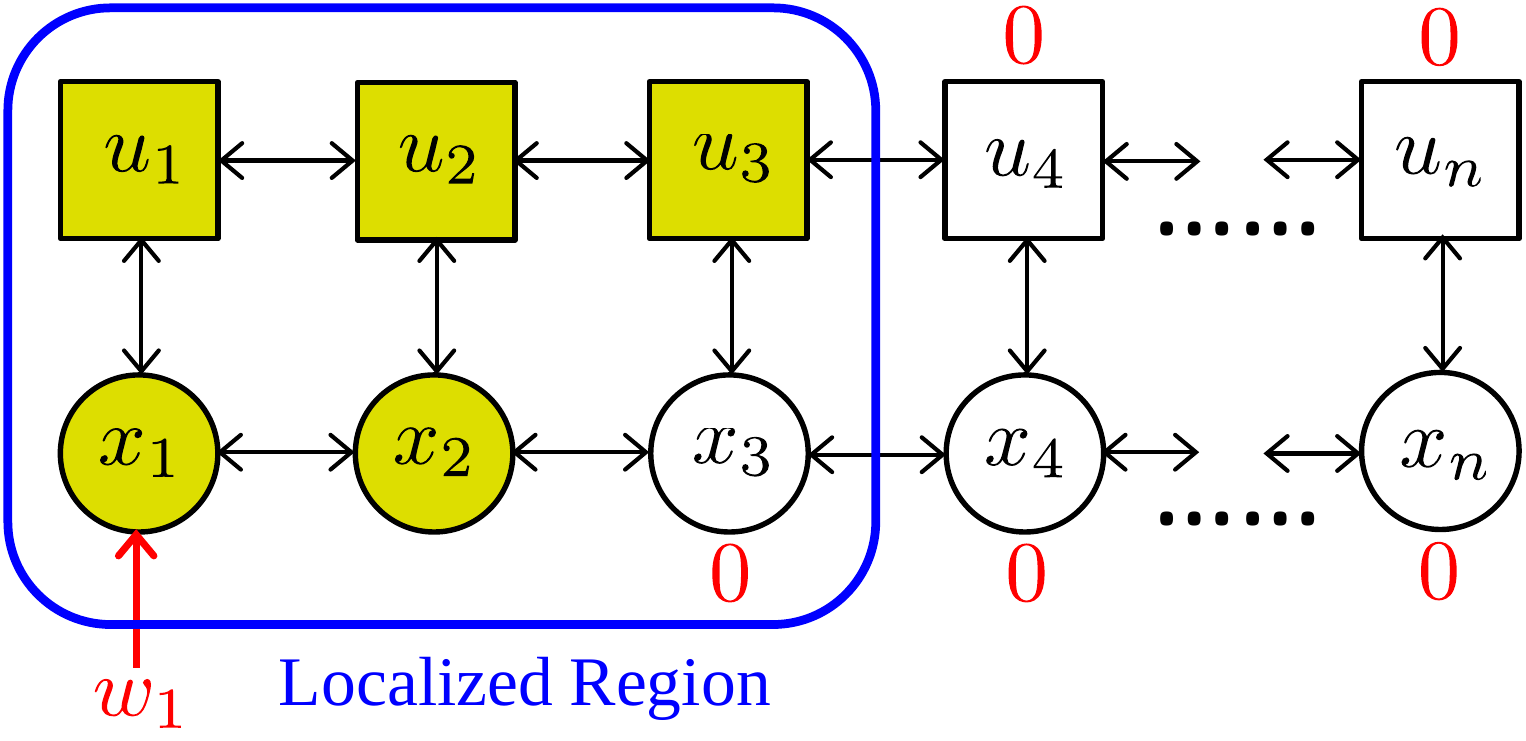}
     \caption{The localized region for $w_1$}
     \label{fig:localized}
\end{figure}
We illustrate the idea of dimension reduction through a simple example below.
\begin{example}
Consider a chain of $n$ LTI subsystems with a single disturbance $\tf w_1$ directly affecting state $\tf x_1$, as shown in Figure \ref{fig:localized}. We analyze the first column of the LLQR problem \eqref{eq:llqr}, which is given by 
\begin{subequations} \label{eq:llqr-r}
\begin{align}
\underset{\{\Sub{\tf R}{:}{1}, \Sub{\tf M}{:}{1}\}}{\text{minimize}} \quad & || \begin{bmatrix} C_1 & D_{12} \end{bmatrix} \begin{bmatrix} \Sub{\tf R}{:}{1} \\ \Sub{\tf M}{:}{1} \end{bmatrix} ||_{\mathcal{H}_2}^2 \label{eq:llqr-r1}\\
\text{subject to} \quad & \begin{bmatrix} zI - A & -B_2 \end{bmatrix} \begin{bmatrix} \Sub{\tf R}{:}{1} \\ \Sub{\tf M}{:}{1} \end{bmatrix} = \Sub{I}{:}{1} \label{eq:llqr-r2} \\
& \begin{bmatrix} \Sub{\tf R}{:}{1} \\ \Sub{\tf M}{:}{1} \end{bmatrix} \in \Sub{\Ell}{:}{1} \cap \FT \cap \frac{1}{z} \mathcal{RH}_\infty. \label{eq:llqr-r3}
\end{align}
\end{subequations}
The optimization variables are the system response from the disturbance $\tf w_1$ to the global state $\tf x$ and global control action $\tf u$.  Assume that the locality constraint imposed on $\tf R$ and $\tf M$ are such that the closed loop response due to $w_1$ satisfies $x_3 = \dots = x_n = 0$ and $u_4 = \dots = u_n = 0$. In this case, we can solve \eqref{eq:llqr-r} using only the information contained within the localized region shown in Figure \ref{fig:localized}. Intuitively, as long as the boundary constraint $x_3 = 0$ is enforced, the perturbation on the system caused by $\tf w_1$ cannot propagate beyond the localized region (this intuition is formalized in Appendix \ref{sec:reduction}). 
Note that the complexity to analyze the system response of $\tf w_1$ is completely independent of the size of the global system.  Even more dramatically, one could replace the section of the system outside of the localized region with any other set of dynamics, and the localized synthesis problem \eqref{eq:llqr-r3} addressing the system response to $\tf w_1$ would remain unchanged.  This localized disturbance-wise analysis can then be performed independently for all other disturbances by exploiting the superposition principle satisfied by LTI systems.
\end{example}

\subsection{Row-wise Separable SLS Problems}
The technique described in the previous subsections can be readily extended to row-wise separable SLS problems -- we briefly introduce the corresponding definitions and end with an example of a row-wise separable SLS problem.

Consider a state estimation SLS problem \cite{2015_PartI} given by
\begin{subequations} \label{eq:sf_lop20}
\begin{align} 
\underset{\ttf{\Phi}}{\text{minimize }} \quad & g(\ttf{\Phi}) \label{eq:lop-10}\\
\text{subject to } \quad & \ttf{\Phi} \tf Z_{AC} = I \label{eq:lop-20} \\
& \ttf{\Phi} \in \s. \label{eq:lop-30}
\end{align}
\end{subequations}
with $\ttf{\Phi} = \begin{bmatrix} \tf R & \tf N \end{bmatrix}$ and $\tf Z_{AC} = \begin{bmatrix} zI-A^\top & -C_2^\top \end{bmatrix}^\top$.
Let $\{ \rowps_1, \dots, \rowps_q\}$ be a partition of the set $\{1, \dots, n_x\}$. The row-wise separability of the SLO and SLC in \eqref{eq:sf_lop20} are defined as follows.
\begin{definition} \label{dfn:dd_row}
The system level objective $g(\ttf{\Phi})$ in \eqref{eq:lop-10} is said to be \emph{row-wise separable} with respect to the row-wise partition $\{ \rowps_1, \dots, \rowps_q\}$ if
\begin{equation}
g(\ttf{\Phi}) = \sum_{j=1}^q g_j(\Sub{\ttf{\Phi}}{\rowps_j}{ : }) \label{eq:decom_obj_row}
\end{equation}
for some functionals $g_j(\cdot)$ for $j = 1, \dots, q$. 
\end{definition}
\begin{definition} \label{dfn:s_row}
The system level constraint $\s$ in \eqref{eq:lop-30} is said to be \emph{row-wise separable} with respect to the row-wise partition $\{ \rowps_1, \dots, \rowps_p\}$ if 
\begin{equation}
\ttf{\Phi} \in \s \quad \text{if and only if} \quad \Sub{\ttf{\Phi}}{ \rowps_j }{:} \in \s_j \text{  for  } j = 1, \dots, p 
\label{eq:gopt-decom_scon}
\end{equation}
for some sets $\s_j$ for $j = 1, \dots, p$.
\end{definition}

An example of a SLS problem satisfying row-wise separabaility is the localized distributed Kalman filter (LDKF) problem in \cite{2015_Wang_LDKF} (naturally seen as the dual to the LLQR problem), which can be posed as
\begin{eqnarray}
\underset{\{\tf R, \tf N\}}{\text{minimize    }} && || \begin{bmatrix} \tf R & \tf N \end{bmatrix} \begin{bmatrix} B_1 \\ D_{21} \end{bmatrix} ||_{\mathcal{H}_2}^2 \nonumber\\
\text{subject to } && \begin{bmatrix} \tf R & \tf N \end{bmatrix} \begin{bmatrix} zI-A \\ -C_2 \end{bmatrix} = I \nonumber \\
&& \begin{bmatrix} \tf R & \tf N \end{bmatrix} \in \begin{bmatrix} \Ell_{R} & \Ell_{N} \end{bmatrix} \cap \FT \cap \frac{1}{z} \mathcal{RH}_\infty. \nonumber
\end{eqnarray}

\section{Partially Separable SLS Problems} \label{sec:of}
In this section, we define the notion of partially-separable objectives and constraints, and show how these properties can be exploited to solve the general form of the localized SLS problem \eqref{eq:localized} in a localized way using distributed optimization techniques such as ADMM.\footnote{Note that \eqref{eq:localized} is neither column-wise nor row-wise separable due to the coupling constraints \eqref{eq:main_1} and \eqref{eq:main_2}.} 
We further show that many constrained optimal control problems are partially separable. Examples include the (localized) $\mathcal{H}_2$ optimal control problem with joint sensor and actuator regularization and the (localized) mixed $\mathcal{H}_2 / \mathcal{L}_1$ optimal control problem. 

Similar to the organization of the previous section, we begin by defining the partial separability of a generic optimization problem. We then specialize our discussion to the partially separable SLS problem \eqref{eq:localized}, before ending this section with a discussion of computational techniques that can be used to accelerate the ADMM algorithm as applied to solving partially separable SLS problems.

\subsection{Partially Separable Problems} \label{sec:gopt-o}
We begin with a generic optimization problem given by
\begin{subequations} \label{eq:gopt-o}
\begin{align} 
\underset{\ttf{\Phi}}{\text{minimize }} \quad & g(\ttf{\Phi}) \label{eq:gopt-o1}\\
\text{subject to } \quad & \ttf{\Phi} \in \s, \label{eq:gopt-o2}
\end{align}
\end{subequations}
where $\ttf{\Phi}$ is a $m \times n$ transfer matrix, $g(\cdot)$ a functional objective, and $\s$ a set constraint. We assume that the optimization problem \eqref{eq:gopt-o} can be written in the form
\begin{subequations} \label{eq:gopt-od}
\begin{align} 
\underset{\ttf{\Phi}}{\text{minimize }} \quad & \rowdec{g}(\ttf{\Phi}) + \coldec{g}(\ttf{\Phi}) \label{eq:gopt-od1}\\
\text{subject to } \quad & \ttf{\Phi} \in \rowdec{\s} \cap \coldec{\s}, \label{eq:gopt-od2}
\end{align}
\end{subequations}
where $\rowdec{g}(\cdot)$ and $\rowdec{\s}$ are row-wise separable and $\coldec{g}(\cdot)$ and $\coldec{\s}$ are column-wise separable. If optimization problem \eqref{eq:gopt-o} can be written in the form of problem \eqref{eq:gopt-od}, then optimization problem \eqref{eq:gopt-o} is called \emph{partially separable}.

Due to the coupling between the row-wise and column-wise separable components, optimization problem \eqref{eq:gopt-od} admits neither a row-wise nor a column-wise decomposition. Our strategy is to instead use distributed optimization techniques such as ADMM to decouple the row-wise separable component and the column-wise separable components: we show that this leads to iterate subproblems that are column/row-wise separable, allowing the results of the previous section to be applied to solve the iterate subproblems using localized and parallel computation. 

\begin{definition} \label{dfn:do}
Let $\{ \rowps_1, \dots, \rowps_{q}\}$ be a partition of the set $\{1, \cdots, m\}$ and $\{ \colps_1, \dots, \colps_{p}\}$ be a partition of the set $\{1, \cdots, n\}$. The functional objective $g(\ttf{\Phi})$ in \eqref{eq:gopt-o1} is said to be \emph{partially separable} with respect to the row-wise partition $\{ \rowps_1, \dots, \rowps_{q}\}$ and the column-wise partition $\{ \colps_1, \dots, \colps_{p}\}$ if $g(\ttf{\Phi})$ can be written as the sum of $\rowdec{g}(\ttf{\Phi})$ and $\coldec{g}(\ttf{\Phi})$, where $\rowdec{g}(\ttf{\Phi})$ is row-wise separable with respect to the row-wise partition $\{ \rowps_1, \dots, \rowps_{q}\}$ and $\coldec{g}(\ttf{\Phi})$ is column-wise separable with respect to the column-wise partition $\{ \colps_1, \dots, \colps_{p} \}$. Specifically, we have 
\begin{eqnarray}
g(\ttf{\Phi}) &=& \rowdec{g}(\ttf{\Phi}) + \coldec{g}(\ttf{\Phi}) \nonumber\\
\rowdec{g}(\ttf{\Phi}) &=& \sum_{j=1}^{q} \rowdec{g}_j(\Sub{\ttf{\Phi}}{\rowps_j}{:}) \nonumber\\
\coldec{g}(\ttf{\Phi}) &=& \sum_{j=1}^{p} \coldec{g}_j(\Sub{\ttf{\Phi}}{:}{\colps_j}) \label{eq:obj_decom2}
\end{eqnarray}
for some functionals $\rowdec{g}_j(\cdot)$ for $j = 1, \dots, q$, and $\coldec{g}_j(\cdot)$ for $j = 1, \dots, p$.  
\end{definition}
\begin{remark}
The column and row-wise separable objective functions defined in Definition \ref{dfn:dd} and \ref{dfn:dd_row}, respectively, are special cases of partially separable objective functions.
\end{remark}

\begin{example}
Let $\Phi$ be the optimization variable, and $E$ and $F$ be two matrices of compatible dimension. The objective function
\begin{equation}
\| E \Phi \|_{\mathcal{F}}^2 + \| \Phi F \|_{\mathcal{F}}^2 \nonumber
\end{equation}
is then partially separable. Specifically, the first term is column-wise separable, and the second term is row-wise separable.
\end{example}

\begin{definition} \label{dfn:ss}
The set constraint $\s$ in \eqref{eq:gopt-o2} is said to be \emph{partially separable} with respect to the row-wise partition $\{ \rowps_1, \dots, \rowps_{q}\}$ and the column-wise partition $\{ \colps_1, \dots, \colps_{p}\}$ if $\s$ can be written as the intersection of two sets $\rowdec{\s}$ and $\coldec{\s}$, where $\rowdec{\s}$ is row-wise separable with respect to the row-wise partition $\{ \rowps_1, \dots, \rowps_{q}\}$ and $\coldec{\s}$ is column-wise separable with respect to the column-wise partition $\{ \colps_1, \dots, \colps_{p} \}$. Specifically, we have 
\begin{eqnarray}
\s &=& \rowdec{\s} \cap \coldec{\s} \nonumber\\
\ttf{\Phi} \in \rowdec{\s} \quad &\Longleftrightarrow& \quad \Sub{\ttf{\Phi}}{\rowps_j}{:} \in \rowdec{\s}_j, \forall j \nonumber\\
\ttf{\Phi} \in \coldec{\s} \quad &\Longleftrightarrow& \quad \Sub{\ttf{\Phi}}{:}{\colps_j} \in \coldec{\s}_j, \forall j \label{eq:sdecom}
\end{eqnarray}
for some sets $\rowdec{\s}_j$ for $j = 1, \dots, q$ and $\coldec{\s}_j$ for $j = 1, \dots, p$.
\end{definition}
\begin{remark}
The column and row-wise separable constraints defined in Definition \ref{dfn:s} and \ref{dfn:s_row}, respectively, are special cases of partially separable constraints.
\end{remark}

\begin{example} [Affine Subspace]
The affine subspace constraint
\begin{equation}
\tf G \ttf{\Phi} = \tf H \quad \text{and} \quad \ttf{\Phi} \tf E = \tf F \nonumber
\end{equation}
is partially separable. Specifically, the first constraint is column-wise separable, and the second is row-wise separable.
\end{example}
\begin{example} [Induced Matrix Norm]
The induced $1$-norm of a matrix $\Phi$ is given by
\begin{equation}
\| \Phi \|_{1} = \underset{ 1 \leq j \leq n}{\text{max}} \,\, \sum_{i=1}^m | \Phi_{ij} |, \nonumber
\end{equation}
and the induced $\infty$-norm of a matrix is given by
\begin{equation}
\| \Phi \|_{\infty} = \underset{ 1 \leq i \leq m}{\text{max}} \,\, \sum_{j=1}^n | \Phi_{ij} |. \nonumber
\end{equation}
The following constraint is partially separable:
\begin{equation}
\| \Phi \|_{1} \leq \gamma_1 \quad \text{and} \quad \| \Phi \|_{\infty} \leq \gamma_{\infty}. \nonumber
\end{equation}
Specifically, the first term is column-wise separable, and the second term is row-wise separable.
\end{example}

Assume now that the objective and the constraints in optimization problem \eqref{eq:gopt-o} are both partially separable with respect to a row-wise partition $\{ \rowps_1, \dots, \rowps_{q}\}$ and a column-wise partition $\{ \colps_1, \dots, \colps_{p}\}$. In this case, the optimization problem \eqref{eq:gopt-o} is partially separable and can be rewritten in the form specified by equation \eqref{eq:gopt-od}. We now propose an ADMM based algorithm that exploits the partially separable structure of the optimization problem.  As we show in the next subsection, when this method is applied to the SLS problem \eqref{eq:localized}, the iterate subproblems are column/row-wise separable, allowing us to apply the methods described in the previous section.  In particular, if locality constraints are imposed, iterate subproblems can be solved using localized and parallel computation.

Let $\ttf{\Psi}$ be a duplicate of the optimization variable $\ttf{\Phi}$. We define the extended-real-value functionals $\rowdec{h}(\ttf{\Phi})$ and $\coldec{h}(\ttf{\Psi})$ by
\begin{eqnarray}
&& \rowdec{h}(\ttf{\Phi}) = \left\{ \begin{array}{ll}
\rowdec{g}(\ttf{\Phi}) & \textrm{if } \ttf{\Phi} \in \rowdec{\s}\\
\infty & \textrm{otherwise}
\end{array} \right. \nonumber\\
&& \coldec{h}(\ttf{\Psi}) = \left\{ \begin{array}{ll}
\coldec{g}(\ttf{\Psi}) & \textrm{if } \ttf{\Psi} \in \coldec{\s}\\
\infty & \textrm{otherwise.}
\end{array} \right. \label{eq:extended}
\end{eqnarray}
Problem \eqref{eq:gopt-od} can then be reformulated as
\begin{eqnarray}
\underset{\{ \ttf{\Phi}, \ttf{\Psi} \}}{\text{minimize    }} && \rowdec{h}(\ttf{\Phi}) + \coldec{h}(\ttf{\Psi}) \nonumber\\
\text{subject to }  && \ttf{\Phi} = \ttf{\Psi}. \label{eq:sub3}
\end{eqnarray}
Problem \eqref{eq:sub3} can be solved using ADMM \cite{BoydADMM}:
\begin{subequations}\label{eq:admm}
\begin{align}
\ttf{\Phi}^{k+1} &= \underset{\ttf{\Phi}}{\text{argmin }} \Big( \rowdec{h}(\ttf{\Phi}) + \frac{\rho}{2} || \ttf{\Phi} - \ttf{\Psi}^{k} + \ttf{\Lambda}^{k}||_{\mathcal{H}_2}^2 \Big) \nonumber \\
\label{eq:admm1} \\
\ttf{\Psi}^{k+1} &= \underset{\ttf{\Psi}}{\text{argmin }} \Big( \coldec{h}(\ttf{\Psi}) + \frac{\rho}{2} || \ttf{\Psi} - \ttf{\Phi}^{k+1} - \ttf{\Lambda}^{k}||_{\mathcal{H}_2}^2 \Big) \nonumber \\
\label{eq:admm2} \\
\ttf{\Lambda}^{k+1} &= \ttf{\Lambda}^{k} + \ttf{\Phi}^{k+1} - \ttf{\Psi}^{k+1} \label{eq:admm3}
\end{align}
\end{subequations}
where the square of the $\mathcal{H}_2$ norm is computed as in \eqref{eq:h2_comp}. 



In Appendix \ref{sec:ADMM}, we provide stopping criterion and prove convergence of the ADMM algorithm \eqref{eq:admm1} - \eqref{eq:admm3} to an optimal solution of \eqref{eq:gopt-od} (or equivalently, \eqref{eq:sub3}) under the following assumptions. 
\begin{assumption}
Problem \eqref{eq:gopt-od} has a feasible solution in the relative interior of the set $\s$. \label{as:1}
\end{assumption}
\begin{assumption} 
The functionals $\rowdec{g}(\cdot)$ and $\coldec{g}(\cdot)$ are closed, proper, and convex. \label{as:2}
\end{assumption}
\begin{assumption} 
The sets $\rowdec{\s}$ and $\coldec{\s}$ are closed and convex. \label{as:3}
\end{assumption}

\subsection{Partially Separable SLS Problems} \label{sec:ps}
We now specialize our discussion to the output feedback SLS problem \eqref{eq:localized}. We use
\begin{equation}
\ttf{\Phi} = \begin{bmatrix} \tf R & \tf N \\ \tf M & \tf L \end{bmatrix} \nonumber
\end{equation}
to represent the system response we want to optimize for. Denote $\tf Z_{AB}$ the transfer matrix $\begin{bmatrix} zI-A & -B_2 \end{bmatrix}$ and $\tf Z_{AC}$ the transfer matrix $\begin{bmatrix} zI-A^\top & -C_2^\top \end{bmatrix}^\top$. Let $J_B$ be the matrix in the right-hand-side of \eqref{eq:main_1} and $J_C$ be the matrix in the right-hand-side of \eqref{eq:main_2}. The localized SLS problem \eqref{eq:localized} can then be written as
\begin{subequations} \label{eq:of_lop}
\begin{align} 
\underset{\ttf{\Phi}}{\text{minimize }} \quad & g(\ttf{\Phi}) \label{eq:of-1}\\
\text{subject to } \quad & \tf Z_{AB} \ttf{\Phi} = J_B \label{eq:of-2} \\
& \ttf{\Phi} \tf Z_{AC} = J_C \label{eq:of-3} \\
& \ttf{\Phi} \in \s \label{eq:of-4}
\end{align}
\end{subequations}
with $\s = \Ell \cap \FT \cap \Sother$. Note that as already discussed, the affine subspace constraints \eqref{eq:of-2} and \eqref{eq:of-3} are partially separable with respect to arbitrary column-wise and row-wise partitions, respectively. Thus the SLS problem \eqref{eq:of_lop} is partially separable if the SLO \eqref{eq:of-1} and the SLC $\ttf \Phi \in \Sother$ are partially separable.

In particular, if $\Sother$ is partially separable, we can express the original SLC $\s$ as an intersection of the sets $\rowdec{\s} = \Ell \cap \FT \cap \rowdec{\Sother}$ and $\coldec{\s} = \Ell \cap \FT \cap \coldec{\Sother}$, where $\rowdec{\Sother}$ is a row-wise separable component of $\Sother$ and $\coldec{\Sother}$ a column-wise separable component of $\Sother$.
Note that the locality constraint and the FIR constraint are included in both $\rowdec{\s}$ and $\coldec{\s}$. 
This is the key point to allow the subroutines \eqref{eq:admm1} - \eqref{eq:admm2} of the ADMM algorithm to be solved using the techniques described in the previous section.

To illustrate how the iterate subproblems of the ADMM algorithm \eqref{eq:admm}, as applied to the SLS problem \eqref{eq:of_lop}, can be solved using localized and parallel computation, we analyze in more detail the iterate subproblem \eqref{eq:admm2}, which is given by
\begin{subequations} \label{eq:admm_lop}
\begin{align} 
\underset{\ttf{\Psi}}{\text{minimize }} \quad & \coldec{g}(\ttf{\Psi}) + \frac{\rho}{2} || \ttf{\Psi} - \ttf{\Phi}^{k+1} - \ttf{\Lambda}^{k}||_{\mathcal{H}_2}^2 \label{eq:admmlop1}\\
\text{subject to } \quad & \ttf{\Psi} \in \coldec{\s} \label{eq:admmlop3}
\end{align}
\end{subequations}
with $\coldec{\s} = \Ell \cap \FT \cap \coldec{\Sother}$.
The $\mathcal{H}_2$ norm regularizer in \eqref{eq:admmlop1} is column-wise separable with respect to arbitrary column-wise partition.
As the objective $\coldec{g}(\cdot)$ and the constraint $\coldec{\s}$ are column-wise separable with respect to a given column-wise partition, we can use the column-wise separation technique described in the previous section to decompose \eqref{eq:admm_lop} into parallel subproblems. 
We can then exploit the locality constraint $\Ell$ and use the technique described in Section \ref{sec:reduction-p} and Appendix \ref{sec:reduction} to reduce the dimension of each subproblem from global to local scale.  Similarly, iterate subproblem \eqref{eq:admm1} can also be solved using parallel and localized computation by exploiting its row-wise separability.  Finally, update equation \eqref{eq:admm3} trivially decomposes element-wise since it is a matrix addition. 

We now give some examples of partially separable SLOs.

\begin{example}[Norm Optimal Control] \label{ex:b1}
Consider the SLO of the distributed optimal control problem in \eqref{eq:trad_obj2}, with the norm given by either the square of the $\mathcal{H}_2$ norm or the element-wise $\ell_1$ norm defined in Example \ref{ex:l1}. Suppose that there exists a permutation matrix $\Pi$ such that the matrix $\begin{bmatrix} B_1^\top & D_{21}^\top \end{bmatrix} \Pi$ is block diagonal. Using a similar argument as in Example \ref{ex:perm}, we can find a column-wise partition to decompose the SLO in a column-wise manner. Similarly, suppose that there exists a permutation matrix $\Pi$ such that the matrix $\begin{bmatrix} C_1 & D_{12} \end{bmatrix} \Pi$ is block diagonal. We can find a row-wise partition to decompose the SLO in a row-wise manner. In both cases, the SLO is column/row-wise separable and thus partially separable.
\end{example}
\begin{example}[Sensor and Actuator Norm]
Consider the weighted actuator norm defined in \cite{MC_CDC14, 2015_Wang_Reg}, which is given by
\begin{equation}
|| \mu \begin{bmatrix} \tf M & \tf L \end{bmatrix} ||_{\mathcal{U}} = \sum_{i=1}^{n_u} \mu_i || e_i^\top \begin{bmatrix} \tf M & \tf L \end{bmatrix} ||_{\mathcal{H}_2} \label{eq:act_norm} 
\end{equation}
where $\mu$ is a diagonal matrix with $\mu_i$ being its $i$th diagonal entry, and $e_i$ is a unit vector with $1$ on its $i$th entry and $0$ elsewhere.The actuator norm \eqref{eq:act_norm} can be viewed as an infinite dimensional analog to the weighted $\ell_1 / \ell_2$ norm, also known as the group lasso \cite{2006_grouplasso} in the statistical learning literature. Adding this norm as a regularizer to a SLS problem induces row-wise sparse structure in the the transfer matrix $\begin{bmatrix} \tf M & \tf L \end{bmatrix}$. Recall from Theorem \ref{thm:1} that the controller achieving the desired system response can be implemented by \eqref{eq:implementation}. If the $i$th row of the transfer matrix $\begin{bmatrix} \tf M & \tf L \end{bmatrix}$ is identically zero, then the $i$th component of the control action $u_i$ is always equal to zero, and therefore the actuator at node $i$ corresponding to control action $u_i$ can be removed without changing the closed loop response. It is clear that the actuator norm defined in \eqref{eq:act_norm} is row-wise separable with respect to arbitrary row-wise partition. This still holds true when the actuator norm is defined by the $\ell_1 / \ell_\infty$ norm. Similarly, consider the weighted sensor norm given by
\begin{equation}
|| \begin{bmatrix} \tf N \\ \tf L \end{bmatrix} \lambda ||_{\mathcal{Y}} = \sum_{i=1}^{n_y} \lambda_i || \begin{bmatrix} \tf N \\ \tf L \end{bmatrix} e_i ||_{\mathcal{H}_2} \label{eq:sen_norm} 
\end{equation}
where $\lambda$ is a diagonal matrix with $\lambda_i$ being its $i$th diagonal entry. The sensor norm \eqref{eq:sen_norm}, when added as a regularizer, induces column-wise sparsity in the transfer matrix $\begin{bmatrix} \tf N^\top & \tf L^\top \end{bmatrix}^\top$. Using the controller implementation \eqref{eq:implementation}, the sensor norm can therefore be viewed as regularizing the number of sensors used by a controller. For instance, if the $i$th column of the transfer matrix $\begin{bmatrix} \tf N^\top & \tf L^\top \end{bmatrix}^\top$ is identically zero, then the sensor at node $i$ and its corresponsing measurement $y_i$ can be removed without changing the closed loop response. The sensor norm defined in \eqref{eq:act_norm} is column-wise separable with respect to any column-wise partition.
\end{example}
\begin{example} [Combination]
From Definition \ref{dfn:do}, it is straightforward to see that the class of partially separable SLOs with respect to the same partitions are closed under summation.
Therefore, we can combine the partially separable SLOs described above, and the resulting SLO is still partially separable. For instance, consider the SLO given by
\begin{eqnarray}
g(\tf R,\tf M,\tf N,\tf L) &=& || \begin{bmatrix} C_1 & D_{12} \end{bmatrix} \begin{bmatrix} \tf R & \tf N \\ \tf M & \tf L \end{bmatrix} \begin{bmatrix} B_1 \\ D_{21} \end{bmatrix} ||_{\mathcal{H}_2}^2 \nonumber\\
&&+ || \mu \begin{bmatrix} \tf M & \tf L \end{bmatrix} ||_{\mathcal{U}} + || \begin{bmatrix} \tf N \\ \tf L \end{bmatrix} \lambda ||_{\mathcal{Y}} \label{eq:h2_sar}
\end{eqnarray}
where $\mu$ and $\lambda$ are the relative penalty between the $\mathcal{H}_2$ performance, actuator and sensor regularizer, respectively. 
If there exists a permutation matrix $\Pi$ such that the matrix $\begin{bmatrix} C_1 & D_{12} \end{bmatrix} \Pi$ is block diagonal, then the SLO \eqref{eq:h2_sar} is partially separable. 
Specifically, the $\mathcal{H}_2$ norm and the actuator regularizer belong to the row-wise separable component, and the sensor regularizer belongs to the column-wise separable component.\footnote{Note that an alternative penalty is proposed in \cite{MC_CDC14} for the design of joint actuation and sensing architectures; it is however more involved to define, and hence we restrict ourself to this alternative penalty for the purposes of illustrating the concept of partially separable SLOs.}
\end{example}

We now provide some examples of partially separable SLCs.
\begin{example} [$\mathcal{L}_1$ Constraint]
The $\mathcal{L}_1$ norm of a transfer matrix is given by its worst case $\ell_\infty$ to $\ell_\infty$ gain.  In particular, the $\mathcal{L}_1$ norm \cite{Munther_L1} of a FIR transfer matrix $\tf G \in \FT$ is given by
\begin{equation}
|| \tf G ||_{\mathcal{L}_1} = \underset{i}{\text{max}} \quad \sum_{j} \sum_{t=0}^T | g_{ij}[t] |.\label{eq:L1_defn}
\end{equation}
We can therefore add the constraint
\begin{equation}
|| \begin{bmatrix} \tf R & \tf N \\ \tf M & \tf L \end{bmatrix} \begin{bmatrix} B_1 \\ D_{21} \end{bmatrix}||_{\mathcal{L}_1} \leq \gamma \label{eq:L1_cons}
\end{equation}
to the optimization problem \eqref{eq:of_lop} for some $\gamma$ to control the worst-case amplification of $\ell_\infty$ bounded signals.  From the definition \eqref{eq:L1_defn}, the SLC \eqref{eq:L1_cons} is row-wise separable with respect to any row-wise partition. 
\end{example}
\begin{example} [Combination]
From Definition \ref{dfn:ss}, the class of partially separable SLCs with respect to the same row and column partitions are closed under intersection, allowing for partially separable SLCS to be combined.
For instance, the combination of a locality constraint $\Ell$, a FIR constraint $\FT$, and an $\mathcal{L}_1$ constraint as in equation \eqref{eq:L1_cons} is partially separable.  This property is extremely useful as it provides a unified framework for dealing with several partially separable constraints at once.
\end{example}

 Using the previously described examples of partially separable SLOs and SLCs, we now consider two partially separable SLS problems: (i) the localized $\mathcal{H}_2$ optimal control problem with joint sensor and actuator regularization, and (ii) the localized mixed $\mathcal{H}_2 / \mathcal{L}_1$ optimal control problem. These two problems are used in Section \ref{sec:simu} as case study examples.

\begin{example}
The localized $\mathcal{H}_2$ optimal control problem with joint sensor and actuator regularization is given by
\begin{subequations} \label{eq:h2SA}
\begin{align} 
\underset{\{\tf R, \tf M, \tf N, \tf L\}}{\text{minimize }} \quad & \eqref{eq:h2_sar} \label{eq:h2SA-1} \\
\text{subject to } \quad & \eqref{eq:main_1} - \eqref{eq:main_stable}  \label{eq:h2SA-2} \\
& \begin{bmatrix} \tf R & \tf N \\ \tf M & \tf L \end{bmatrix} \in \C \cap \Ell \cap \FT. \label{eq:h2SA-3}
\end{align}
\end{subequations}
where $\C$ encodes the information sharing constraints of the distributed controller. If there exists a permutation matrix $\Pi$ such that the matrix $\begin{bmatrix} C_1 & D_{12} \end{bmatrix} \Pi$ is block diagonal, then \eqref{eq:h2SA} is partially separable. 
\end{example}
\begin{remark}
When the penalty of the sensor and actuator norms are zero, problem \eqref{eq:h2SA} reduces to the localized LQG optimal control problem defined and solved in \cite{2015_Wang_H2}.  Further, if the locality and FIR constraints are removed, we recover the standard LQG optimal control problem, which is also seen to be partially separable.
\end{remark}


Next we consider the localized mixed $\mathcal{H}_2/\mathcal{L}_1$ optimal control problem.
\begin{example}
The localized mixed $\mathcal{H}_2/\mathcal{L}_1$ optimal control problem is given by 
\begin{subequations} \label{eq:h2l1}
\begin{align} 
\underset{\{\tf R, \tf M, \tf N, \tf L\}}{\text{minimize }} \quad & || \begin{bmatrix} \tf R & \tf N \\ \tf M & \tf L \end{bmatrix} \begin{bmatrix} B_1 \\ D_{21} \end{bmatrix}||_{\mathcal{H}_2}^2 \label{eq:h2l1-1} \\
\text{subject to } \quad & \eqref{eq:main_1} - \eqref{eq:main_stable}, \eqref{eq:L1_cons}, \eqref{eq:h2SA-3} \label{eq:h2l1-2}
\end{align}
\end{subequations}
which is partially separable.
\end{example}

The localized mixed $\mathcal{H}_2/\mathcal{L}_1$ optimal control problem can be used to design the tradeoff between average and worst-case performance, as measured by the $\mathcal{H}_2$ and $\mathcal{L}_1$ norms of the closed loop system, respectively. 


\subsection{Analytic Solution and Acceleration}
Suppose that the Assumptions \ref{as:1} - \ref{as:3} in Section \ref{sec:ps} hold. The ADMM algorithm presented in \eqref{eq:admm} is a special case of the proximal algorithm \cite{proximal, proximalBoyd, BoydADMM}. For certain type of objective functionals $\rowdec{h}(\cdot)$ and $\coldec{h}(\cdot)$, the proximal operators can be evaluated analytically (see Ch. 6 of \cite{proximalBoyd}). In this situation, we only need to evaluate the proximal operators \emph{once}, and iterate \eqref{eq:admm1} and \eqref{eq:admm2} in closed form. This improves the overall computation time significantly. 
We explain how to express the solutions of \eqref{eq:admm1} and \eqref{eq:admm2} using proximal operators in detail in Appendix \ref{sec:appendix}. Here we list a few examples for which the proximal operators can be evaluated analytically. 

\begin{example} \label{ex:analytic1}
Consider the LLQG problem proposed in \cite{2015_Wang_H2}.
When the global optimization problem is decomposed into parallel subproblems, each subproblem is a convex quadratic program restricted to an affine set. In this case, the proximal operator is an affine function \cite{2015_Wang_H2, proximalBoyd, BoydADMM}. We only need to calculate this affine function once. The iteration in \eqref{eq:admm1} - \eqref{eq:admm3} can then be carried out using multiple matrix multiplications in the reduced dimension, which significantly improves the overall computation time.
\end{example}
\begin{example} \label{ex:analytic2}
Consider the LLQR problem with actuator regularization \cite{2015_Wang_Reg}, which is the state feedback version of \eqref{eq:h2SA}. The column-wise separable component of this problem is identical to that of the LLQG example, and therefore the update subproblem \eqref{eq:admm2} can be solved using matrix multiplication, as described in the previous example.  We further showed in \cite{2015_Wang_Reg} that the row-wise separable component of this problem can be simplified to multiple unconstrained optimization problems, each with proximal operators given by vectorial soft-thresholding \cite{proximalBoyd}.  This offers an efficient way to solve the update supbroblem \eqref{eq:admm1}.
\end{example}

We end this section by noting that although our focus has been on how ADMM can be used to efficiently solve the localized SLS problem \eqref{eq:localized},  there exist other distributed algorithms that can exploit partial separability to improve computational efficiency.  For instance, if either $\rowdec{g}(\cdot)$ or $\coldec{g}(\cdot)$ is strongly convex, we can use the alternating minimization algorithm (AMA) \cite{AMA} to simplify the ADMM algorithm. 




\section{Simulations} \label{sec:simu}
In this section, we apply the localized $\mathcal{H}_2$ optimal control with sensor actuator regularization \eqref{eq:h2SA} problem and the localized mixed $\mathcal{H}_2 / \mathcal{L}_1$ optimal control problem \eqref{eq:h2l1} to a power system inspired example.
After introducing the power system model, we show that the localized $\mathcal{H}_2$ controller, with its additional locality, FIR and communication delay constraints can achieve comparable closed loop performance to that of a centralized $\mathcal{H}_2$ optimal controller.  We further demonstrate the scalability of the proposed method by synthesizing a localized $\mathcal{H}_2$ controller for a randomized heterogeneous networked system with $12800$ states using a single laptop; for such a large system, neither the centralized nor the distributed optimal controllers can be computed. We then solve the localized $\mathcal{H}_2$ with joint sensor and actuator regularization to co-design an output feedback controller and its actuation/sensing architecture. Finally, we solve the localized mixed $\mathcal{H}_2 / \mathcal{L}_1$ optimal control problem, thus identifying the achievable tradeoff curve between the average and worst-case performance of the closed loop system.

\subsection{Power System Model}
\begin{figure}[ht!]
      \centering
      \subfigure[Interconnected topology]{%
      \includegraphics[width=0.36\columnwidth]{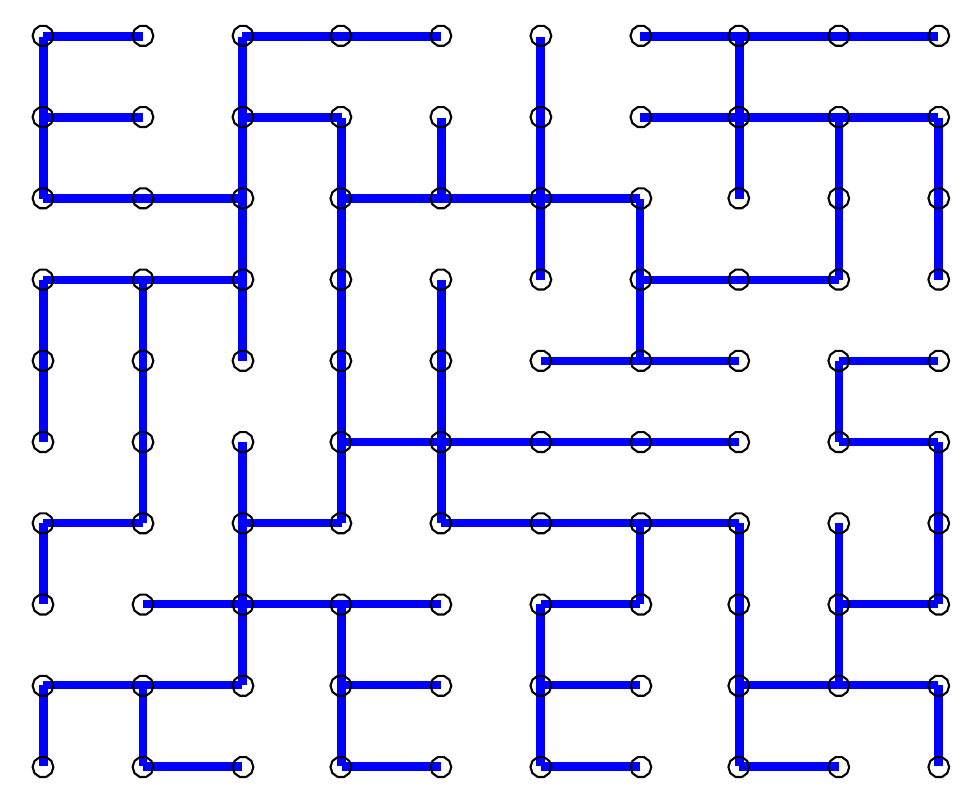}
      \label{fig:mesh}}      
      \subfigure[Interaction between neighboring subsystems]{%
      \includegraphics[width=0.45\columnwidth]{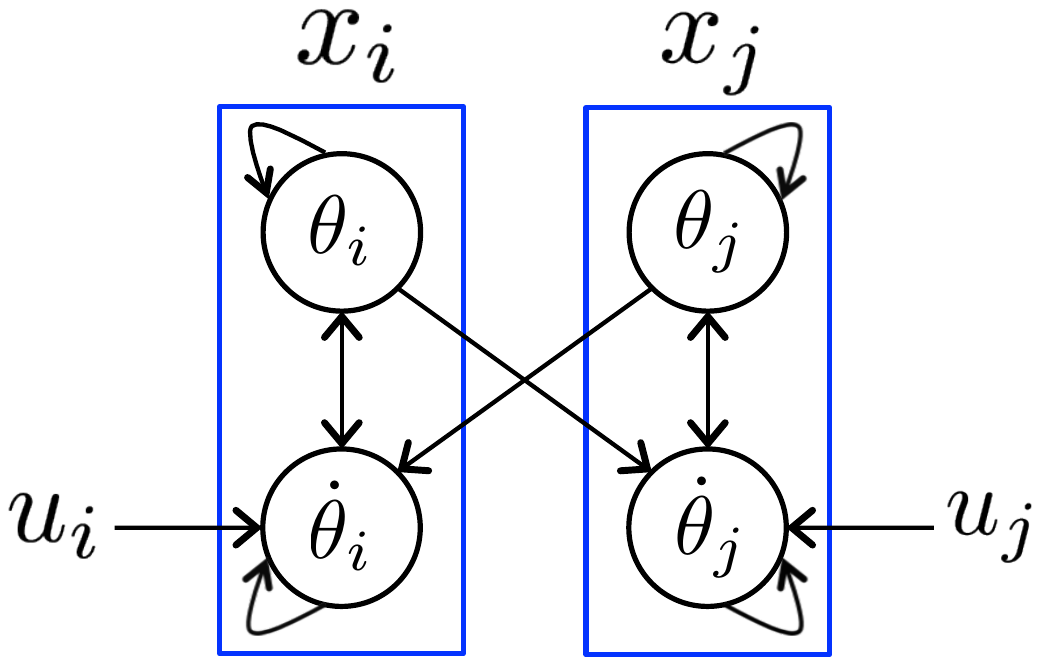}
      \label{fig:inter}} 
      \caption{Simulation example interaction graph.}
\end{figure}

We begin with a randomized spanning tree embedded on a $10 \times 10$ mesh network representing the interconnection between subsystems. The resulting interconnected topology is shown in Fig. \ref{fig:mesh} -- we assume that all edges are undirected.
The dynamics of each subsystem is given by the discretized swing equations. Consider the swing dynamics 
\begin{equation}
m_i \ddot{\theta_i} + d_i \dot{\theta_i} = -\sum_{j \in \mathcal{N}_i} k_{ij} (\theta_i - \theta_j) + w_i + u_i \label{eq:swing}
\end{equation}
where $\theta_i$, $\dot{\theta_i}$, $m_i$, $d_i$, $w_i$, $u_i$ are the phase angle deviation, frequency deviation, inertia, damping, external disturbance, and control action of the controllable load of bus $i$. The coefficient $k_{ij}$ is the coupling term between buses $i$ and $j$. We let $x_i := [\theta_i \quad \dot{\theta_i}]^\top$ be the state of bus $i$ and use $e^{A \Delta t} \approx I + A \Delta t$ to discretize the swing dynamics. Equation \eqref{eq:swing} can then be expressed in the form of \eqref{eq:interconnected} with
\begin{equation}
A_{ii} = \begin{bmatrix} 1 & \quad \Delta{t} \\ -\frac{k_i}{m_i} \Delta{t} & \quad 1-\frac{d_i}{m_i} \Delta{t} \end{bmatrix}, \quad A_{ij} = \begin{bmatrix} 0 & 0 \\ \frac{k_{ij}}{m_i} \Delta{t} & 0 \end{bmatrix}, \nonumber
\end{equation}
and $B_{ii} = \begin{bmatrix} 0 & 1 \end{bmatrix}^\top$. We set $\Delta{t} = 0.2$ and $k_i = \sum_{j \in \mathcal{N}_i} k_{ij}$. In addition, the parameters $k_{ij}, d_i$, and $m_i^{-1}$ are randomly generated and uniformly distributed between $[0.5, 1]$, $[1, 1.5]$, and $[0, 2]$, respectively. The $A$ matrix is then normalized to be marginally stable, i.e., such that the spectral radius of the matrix $A$ is $1$. The interaction between neighboring subsystems of the discretized model is illustrated in Figure \ref{fig:inter}. We initially assume that each subsystem in the power network has a phase measurement unit (PMU), a frequency sensor, and a controllable load that generates $u_i$ -- later in this section we explore the tradeoff between actuation/sensing density and closed loop performance.

From \eqref{eq:swing}, the external disturbance $w_i$ only directly affects the frequency deviation $\dot{\theta_i}$. 
To make the objective functional strongly convex, we introduce a small artificial disturbance on the phase deviation $\theta_i$ as well. 
We assume that the process noise on frequency and phase are uncorrelated AWGNs with covariance matrices given by $I$ and $10^{-4}I$, respectively. In addition, we assume that both the phase deviation and the frequency deviation are measured with some sensor noise. The sensor noise of phase and frequency measurements are uncorrelated AWGNs with covariance matrix given by $10^{-2}I$. We choose equal penalty on the state deviation and control effort, i.e., $\begin{bmatrix} C_1 & D_{12} \end{bmatrix} = I$. 

Based on the above setting, we formulate a $\mathcal{H}_2$ optimal control problem that minimizes the $\mathcal{H}_2$ norm of the transfer matrix from the process and sensor noises to the regulated output. The $\mathcal{H}_2$ norm of the closed loop is given by $13.3169$ when a proper centralized $\mathcal{H}_2$ optimal controller is applied, and $16.5441$ when a strictly proper centralized $\mathcal{H}_2$ optimal controller is applied. In the rest of this section, we normalized the $\mathcal{H}_2$ norm with respect to the proper centralized $\mathcal{H}_2$ optimal controller.



\subsection{Localized $\mathcal{H}_2$ Optimal Control}
The underlying assumption of the centralized optimal control scheme is that the measurement can be transmitted \emph{instantaneously} with \emph{every subsystem} in the network. 
To incorporate realistic communication delay constraint and facilitate the scalability of controller design, we impose additional communication delay constraint, locality constraint, and FIR constraint on the system response. 

For the communication delay constraint $\C$, we assume that each subsystem takes one time step to transmit the information to its neighboring subsystems.  Therefore, if subsystems $i$ and $j$ are $k$-hops apart (as defined by the interaction graph illustrated in Figure \ref{fig:mesh}), then the control action $u_i[t]$ at subsystem $i$ can only use the measurements $y_j[\tau]$ and internal controller state $\beta_j[\tau]$ of subsystem $j$ if $\tau \leq t-k$.

The interaction between subsystems illustrated in Fig. \ref{fig:inter} implies that it takes two discrete time steps for a disturbance at subsystem $j$ to propagate to its neighboring subsystems, and hence the communication speed between sub-controllers is twice as fast as the propagation speed of disturbances through the plant. For the given communication delay constraint $\C$, we use the method in \cite{2015_Wang_Reg} to design the tightest feasible locality constraint $\Ell$. In this example, we can localize the joint effect of the process and sensor noise at subsystem $j$ to a region defined by its two-hop neighbors (where one hop is as defined in terms of the interaction graph of the system). This implies that the sub-controller at node $j$ only needs to transmit its measurements $y_j$ and controller states $\beta_j$ within this localized region, and further only a restricted plant model (as defined by this localized region) is needed to synthesize its corresponding local control policy. 

\begin{figure}[ht!]
      \centering
      \includegraphics[width=0.4\textwidth]{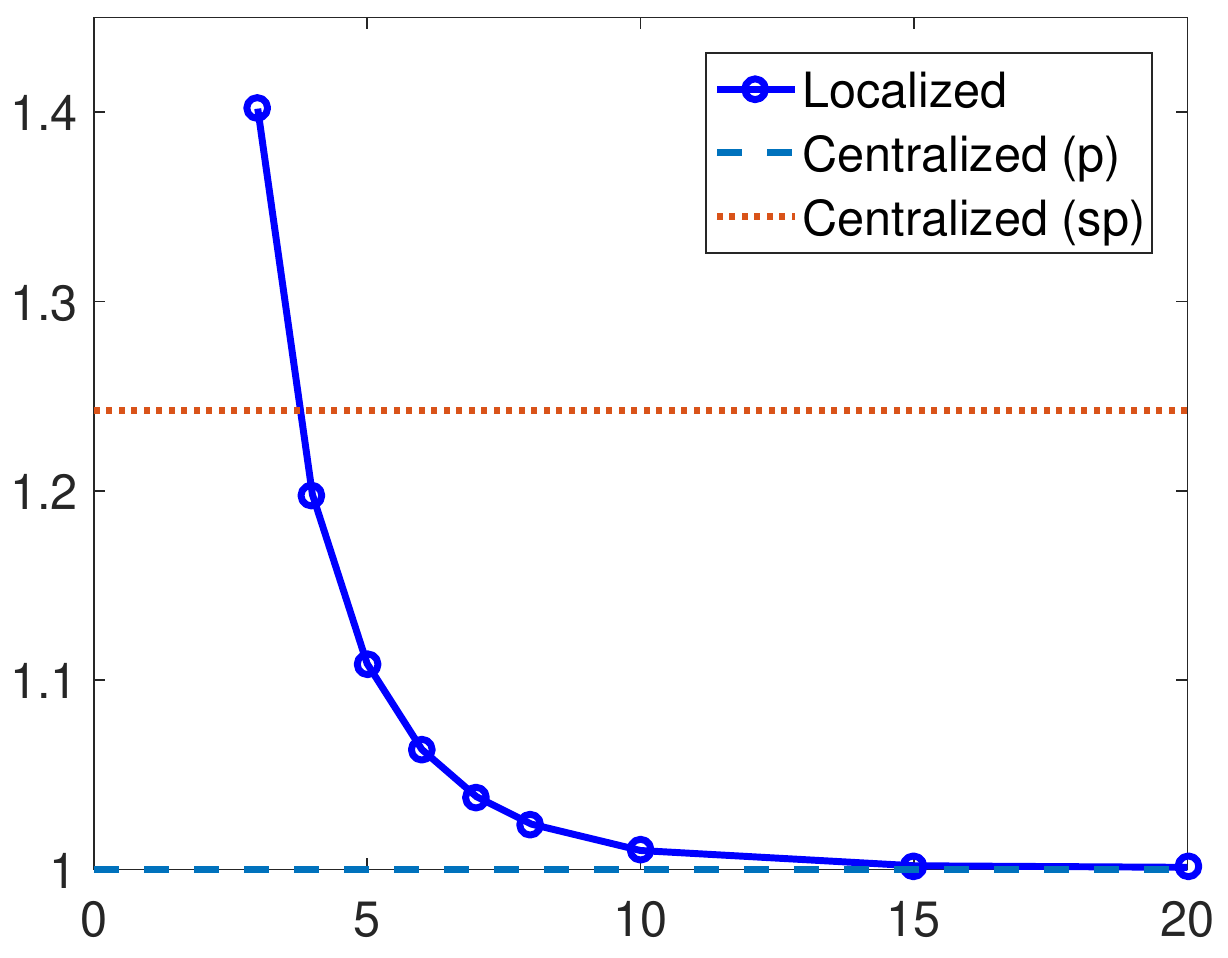}
      \caption{The vertical axis is the normalized $\mathcal{H}_2$ norm of the closed loop when the localized $\mathcal{H}_2$ controller is applied. The localized $\mathcal{H}_2$ controller is subject to the constraint $\C \cap \Ell \cap \FT$. The horizontal axis is the horizon $T$ of the FIR constraint $\FT$, which is also the settling time of an impulse disturbance. We plot the normalized $\mathcal{H}_2$ norm for the centralized unconstrained optimal controller (proper and strictly proper) in the same figure.}
      \label{fig:T}
\end{figure}

Assuming a fixed communication delay constraint $\C$ and locality constraint $\Ell$, we first explore tradeoff between the length $T$ of the FIR constraint $\FT$ and the transient performance. Figure \ref{fig:T} shows the tradeoff curve between the transient performance of the localized $\mathcal{H}_2$ controller and the length $T$ of the FIR constraint.
For the given communication delay constraint $\C$ and the locality constraint $\Ell$, the localized $\mathcal{H}_2$ controller is feasible for FIR constraints $\FT$ whenever $T \geq 3$. When the length of the FIR constraint increases, the $\mathcal{H}_2$ norm of the closed loop converges quickly to the unconstrained optimal value. For instance, for FIR lengths of $T=7,\, 10,$ and $20$, the performance degradation with respect to the unconstrained $\mathcal{H}_2$ optimal controller are given by $3.8\%$, $1.0\%$, and $0.1\%$, respectively. This further means that the performance degradation due to the additional communication delay constraint $\C$ and the locality constraint $\Ell$ is less than $0.1\%$. From Figure \ref{fig:T}, we see that the localized $\mathcal{H}_2$ controller, with its additional communication delay, locality and FIR constraints can achieve similar transient performance to that of an unconstrained centralized (unimplementable) optimal $\mathcal{H}_2$ controller.

To further illustrate the advantages of the localized control scheme, we choose $T = 20$ and compare the localized optimal controller, distributed optimal controller, and the centralized optimal controller in terms of the closed loop performance, the complexity of controller synthesis, and the complexity of controller implementation in Table \ref{Table:2}. The distributed optimal controller is computed using the method described in \cite{2014_Lamperski_H2_journal}, in which we assume the same communication constraint $\C$ as the localized optimal controller. It can be seen that the localized controller is vastly preferable in all aspects, except for a slight degradation in the closed-loop performance.

\begin{table}[h]
 \caption{Comparison Between Centralized, Distributed, and Localized Control}
 \label{Table:2}
\begin{center}
    \begin{tabular}{| l | l | l | l | l |}
    \hline
     & & Cent. & Dist. & Local. \\ \hline
     & Affected region & Global & Global & $2$-hop \\ \cline{2-5}
     Closed Loop & Affected time & Long & Long & $20$ steps \\ \cline{2-5}
     & Normalized $\mathcal{H}_2$ & $1$ & $1.001$ & $1.001$ \\ \hline \hline
     & Complexity & $O(n^3)$ & $\geq O(n^3)$ & $O(n)$ \\ \cline{2-5}
     Synthesis & Plant model & Global & Global & $2$-hop \\ \cline{2-5}
     & Redesign & Offline & Offline & Real-time \\ \hline \hline
     Implement. & Comm. Speed & $\infty$ & $2$ & $2$ \\ \cline{2-5}
     & Comm. Range & Global & Global & $2$-hop \\ \hline 
     \end{tabular}
\end{center}
\end{table}

\begin{figure}[ht!]
      \centering
      \includegraphics[width=0.45\textwidth]{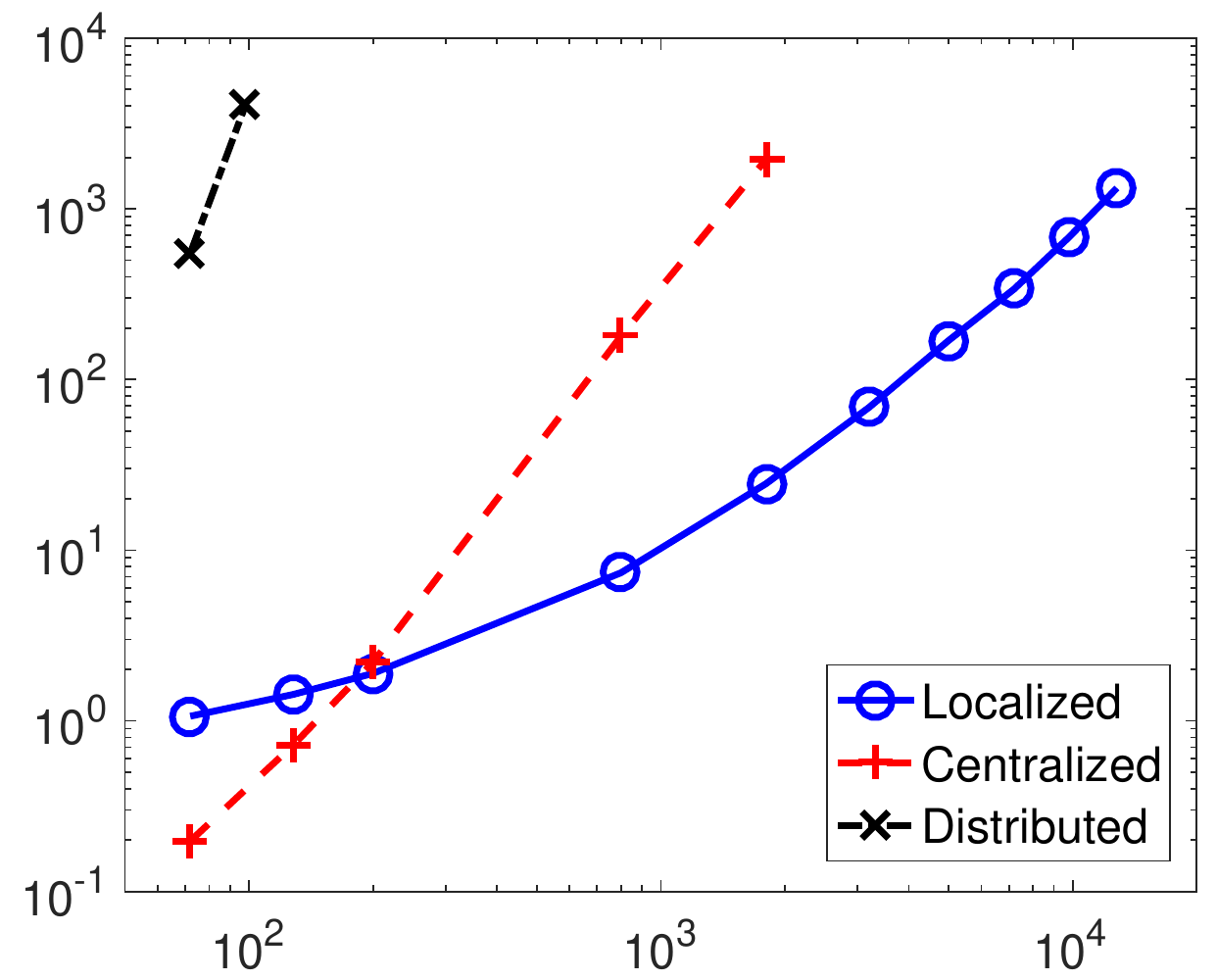}
      \caption{The horizontal axis denotes the number of states of the system, and the vertical axis is the computation time in seconds.}
      \label{fig:Time}
\end{figure}

We now allow the size of the problem to vary and compare the computation time needed to synthesize a centralized, distributed, and localized $\mathcal{H}_2$ optimal controller. We choose $T = 7$ for the localized controller. The empirical relationship obtained between computation time and problem size for different control schemes is illustrated in Figure \ref{fig:Time}. As can be seen in Figure \ref{fig:Time}, the computation time needed for the distributed controller grows rapidly when the size of problem increases. The slope of the line describing the computation time of the centralized controller in the log-log plot of Figure \ref{fig:Time} is $3$, which matches the theoretical complexity $O(n^3)$. The slope for the localized $\mathcal{H}_2$ controller is about $1.4$, which is larger than the theoretical value of $1$. We believe this overhead may caused by other computational issue such as memory management. We note that the computational bottleneck faced in computing our large-scale example arises from using a single laptop to compute the controller (and hence the localized subproblems were solved in serial) -- in practice, if each local subsystem is capable of solving its corresponding localized subproblem, our approach scales to systems of arbitrary size as all computations can be done in parallel. For the largest example that we have, we can compute the optimal localized $\mathcal{H}_2$ controller for a system with $12800$ states in $22$ minutes using a standard laptop. If the computation is parallelized across all $6400$ sub-systems, the synthesis algorithm can be completed within $0.2$ seconds. In contrast, the theoretical time to compute the centralized $\mathcal{H}_2$ optimal controller for the same example is more than a week.

\subsection{Localized $\mathcal{H}_2$ Optimal Control with Joint Sensor and Actuator Regularization}
We now move back to the $10 \times 10$ mesh example shown in Figure \ref{fig:mesh}. In the previous subsection, we assume that each subsystem in the power network had a phase sensor, a frequency sensor, and a controllable load. In practice, the installation of these sensors and actuators is expensive, and we instead like to explore the tradeoff between the closed loop performance of the system and the number of sensors and actuators being used in a principled and tractable manner. A challenging problem is to determine the optimal locations of these sensors and actuators due to its combinatorial complexity. In this subsection, we apply the regularization for design (RFD) \cite{MC_CDC14} framework to jointly design the localized optimal controller and the optimal locations of sensors and actuators in the power network. This is achieved by solving the localized $\mathcal{H}_2$ optimal control problem with joint sensor and actuator regularization \eqref{eq:h2SA}.

In order to allow more flexibility on sensor actuator placement, we increase the localized region of each process and sensor noise from its two-hop neighbors to its four-hop neighbors. This implies that each subsystem $j$ needs to exchange the information up to its four-hop neighbors, and use the restricted plant model within its four-hop neighbors to synthesize the localized $\mathcal{H}_2$ controller. The length of the FIR constraint $\FT$ is increased to $T = 30$. The initial localized $\mathcal{H}_2$ cost is given by $13.3210$, which is $0.03\%$ degradation compared to the idealized centralized $\mathcal{H}_2$ controller. We assume that the relative price between each frequency sensor, PMU, and controllable load are $1$, $100$, and $300$, respectively. This is to model the fact that actuators are typically more expensive than sensors, and that PMUs are typically more expensive than frequency sensors. The price for the same type of sensors and actuators at different locations remains constant.

\begin{figure}[ht!]
      \centering
      \includegraphics[width=0.25\textwidth]{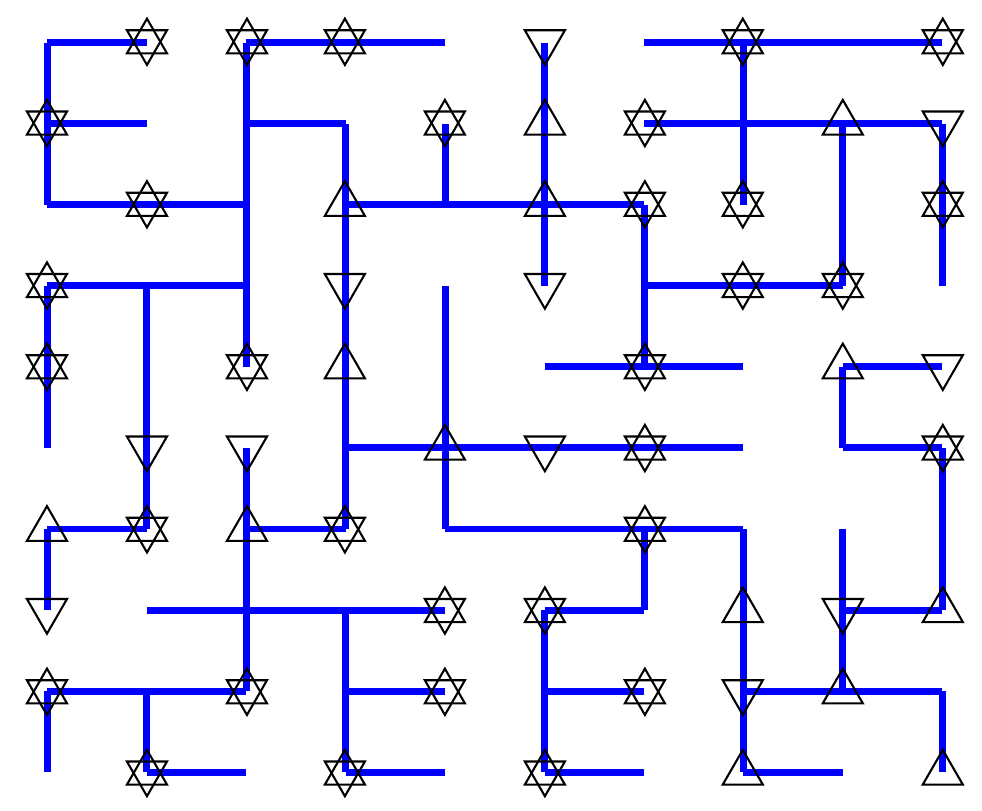}
      \caption{The upward-pointing triangles represent the subsystems in which the PMU is removed. The downward-pointing triangles represent the subsystems in which the controllable load (actuator) is removed.}
      \label{fig:SA}
\end{figure}

We begin with a dense controller architecture composed of $100$ controllable loads, $100$ PMUs, and $100$ frequency sensors, i.e., one of each type of sensor and actuator at each node.  Using optimization problem \eqref{eq:h2SA} to identify the sparsest possible architecture that still satisfies the locality, FIR and communication delay constraints, we are able to remove $43$ controllable loads and $46$ PMUs from the system (no frequency sensors were removed due to the chosen relative pricing). The locations of the removed sensors and actuators are shown in Figure \ref{fig:SA}. We argue that this sensing and actuation interface is very sparse. In particular, we only use $57$ controllable loads to control process noise from $200$ states and sensor noise from $154$ states, while ensuring that the system response to {all} process and sensor disturbances remains both localized and FIR.

The localized $\mathcal{H}_2$ cost for the system with reduced number of sensors and actuators is given by $17.8620$. In comparison, the cost achieved by a proper centralized optimal controller is $16.2280$, and the cost achieved by a strictly proper centralized optimal controller is $18.4707$. Note that as the sensing and actuation interface becomes sparser, the performance gap between the centralized and the localized controller becomes larger. 
Nevertheless, we note that the performance degradation is only $10\%$ compared to the proper centralized optimal scheme implemented using the same sparse controller architecture.

\subsection{Localized Mixed $\mathcal{H}_2/\mathcal{L}_1$ Optimal Control}
\begin{figure}[ht!]
      \centering
      \includegraphics[width=0.35\textwidth]{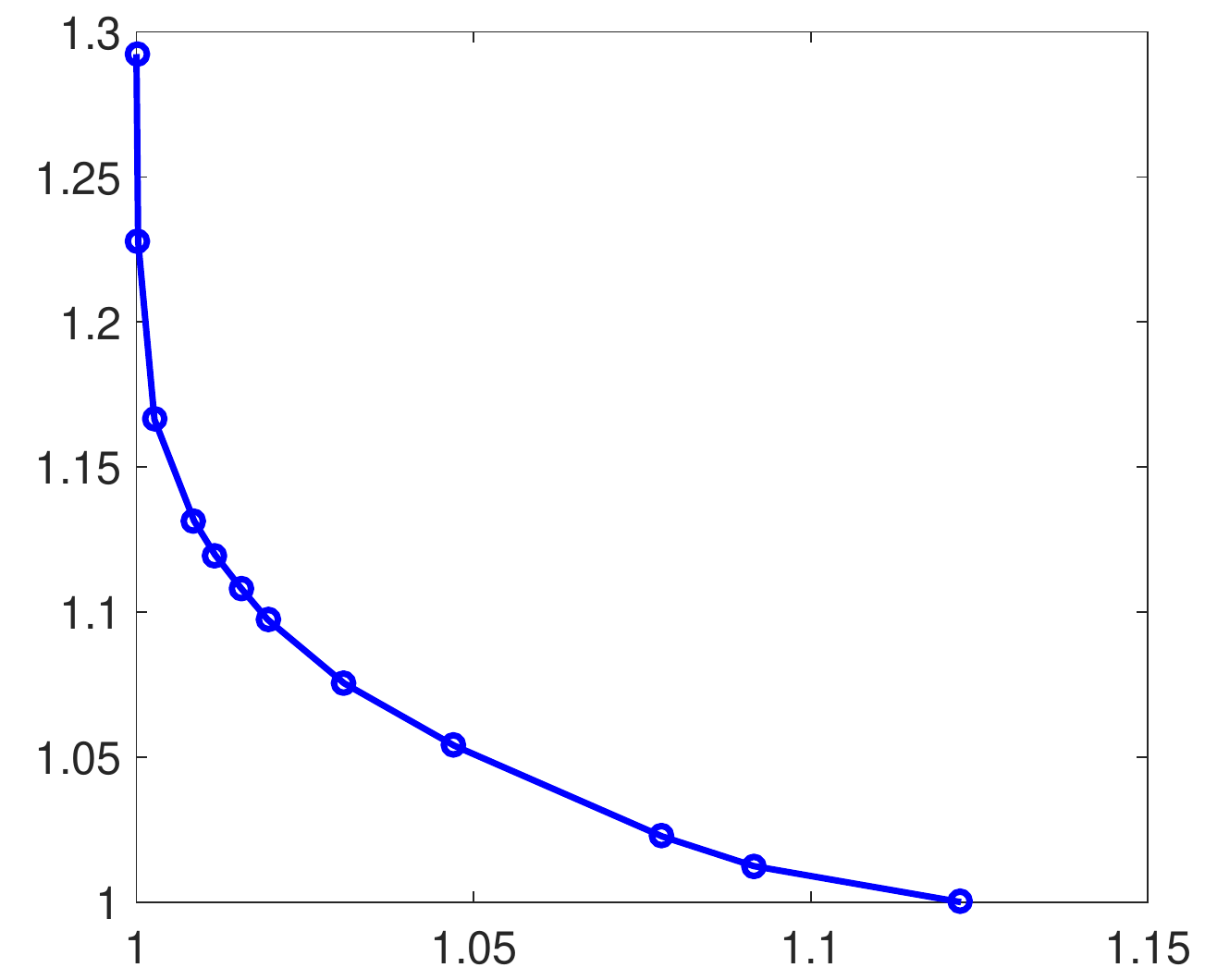}
      \caption{The vertical axis represents the normalized $\mathcal{L}_1$ norm of the closed loop, and the horizontal axis represents the normalized $\mathcal{H}_2$ norm of the closed loop.}
      \label{fig:H2L1}
\end{figure}

Finally, we solve the localized mixed $\mathcal{H}_2/\mathcal{L}_1$ optimal control problem in \eqref{eq:h2l1} on the $10 \times 10$ mesh example shown in Figure \ref{fig:mesh}. We progressively reduce the allowable $\mathcal{L}_1$ gain $\gamma$, as shown in equation \eqref{eq:h2l1}, to explore the tradeoff between average and worst-case performance of a system. We plot the normalized $\mathcal{H}_2$ norm and the normalized $\mathcal{L}_1$ norm in Figure \ref{fig:H2L1}.\footnote{We normalize the $\mathcal{H}_2$ performance with respect to that achieved by the optimal localized $\mathcal{H}_2$ controller, and likewise normalize the $\mathcal{L}_1$ performance with respect to that achieved by the optimal localized $\mathcal{L}_1$ controller.} The left-top point in Figure \ref{fig:H2L1} is the localized $\mathcal{H}_2$ solution. When we start reducing the $\mathcal{L}_1$ sublevel set, the $\mathcal{H}_2$ norm of the closed loop response gradually increases, thus tracing out a tradeoff curve. 


\section{Conclusion} \label{sec:conclusion}
In this paper, we proposed a localized and scalable algorithm to solve a class of constrained optimal control problems for large-scale systems using the system level synthesis (SLS) framework \cite{2015_PartI}. Specifically, we defined suitable notions of separability for SLS problems, and showed that this allowed for them to be solved via parallel (and local) computation.
We then argued that the resulting synthesis algorithms can scale to arbitrarily large systems if suitable parallel computation is available.  We further showed that many constrained optimal control problems of interest satisfy these separability conditions;  these include for instance the (localized) mixed $\mathcal{H}_2 / \mathcal{L}_1$ optimal control problem and the (localized) $\mathcal{H}_2$ optimal control problem with joint sensor and actuator regularization. We further demonstrated the scalability of our algorithm by synthesizing a localized $\mathcal{H}_2$ optimal controller on a randomized heterogeneous power system with $12800$ states, a system for which neither the centralized nort the distributed optimal controllers can be computed in a reasonable amount of time. 
Finally, through the use of several simulation examples, we demonstrated that the localized SLS framework provides a useful tool to explore the various tradeoffs that arise in designing large-scale distributed controllers. 


%

\appendices
\section{Dimension Reduction Algorithm} \label{sec:reduction}

Consider problem \eqref{eq:decom} for a specific $j$. Recall from equation \eqref{eq:desired} that the number of rows of the transfer matrix $\Sub{\ttf{\Phi}}{:}{\colps_j}$ is given by $(n_x + n_u)$. Let $\bar{\rone}_j$ be the largest subset of $\{1, \dots, n_x + n_u\}$ such that the locality constraint $\Sub{\Ell}{\bar{\rone}_j}{\colps_j}$ in \eqref{eq:decom-3} is exactly a zero matrix. When the locality constraint is imposed, we must have $\Sub{\ttf{\Phi}}{\bar{\rone}_j}{\colps_j} = 0$. This component of $\ttf{\Phi}$ can be eliminated from the objective function \eqref{eq:decom-1} and the constraints \eqref{eq:decom-2} - \eqref{eq:decom-3}. Let $\rone_j = \{1, \dots, n_x + n_u\} - \bar{\rone}_j$ be the complement set of $\bar{\rone}_j$. When the above sparsity constraints are imposed on the optimization variable, problem \eqref{eq:decom} can be simplified to
\begin{subequations} \label{eq:decom2}
\begin{align} 
\underset{\Sub{\ttf{\Phi}}{\rone_j}{\colps_j}}{\text{minimize }} \quad & \bar{g}_j(\Sub{\ttf{\Phi}}{\rone_j}{\colps_j}) \label{eq:decom2-1}\\
\text{subject to } \quad & \Sub{\tf Z_{AB}}{:}{\rone_j} \Sub{\ttf{\Phi}}{\rone_j}{\colps_j} = \Sub{I}{ : }{\colps_j} \label{eq:decom2-2} \\
& \Sub{\ttf{\Phi}}{\rone_j}{\colps_j} \in \Sub{\Ell}{\rone_j}{\colps_j} \cap \FT \cap \bar{\Sother}_j \label{eq:decom2-3}
\end{align}
\end{subequations}
where $\bar{g}_j(\cdot)$ and $\bar{\Sother}_j$ are the restrictions of the column-wise objective function $g_j(\cdot)$ and constraint $\Sother_j$, respectively, to the subspace satisfying $\Sub{\ttf{\Phi}}{\bar{\rone}_j}{\colps_j} = 0$. From the definition of the set $\rone_j$, it is straightforward to see that \eqref{eq:decom2} is equivalent to $\eqref{eq:decom}$.

Next, if the system matrices $(A,B_2)$ are sparse, then the transfer matrix $\tf Z_{AB}$ is also sparse. Recall that the number of rows of $\tf Z_{AB}$ is given by $n_x$. Let $\bar{\rtwo}_j$ be the largest subset of $\{1, \dots, n_x\}$ such that the augmented matrix $\begin{bmatrix}\Sub{\tf Z_{AB}}{\bar{\rtwo}_j}{\rone_j} & \Sub{I}{\bar{\rtwo}_j}{\colps_j}\end{bmatrix}$ is a zero matrix. Let $\rtwo_j = \{1, \dots, n_x\} - \bar{\rtwo}_j$ be the complement set of $\bar{\rtwo}_j$. We can further reduce optimization problem \eqref{eq:decom2} to optimization problem \eqref{eq:decom40}.
From the definition of the sets $\rtwo_j$, it follows that \eqref{eq:decom40} is equivalent to \eqref{eq:decom2} and therefore equivalent to \eqref{eq:decom}.


\begin{remark} \label{rmk:1}
The dimension reduction algorithm proposed here is a generalization of the algorithms proposed in our prior works \cite{2014_Wang_CDC, 2015_Wang_H2}. Specifically, the algorithms in \cite{2014_Wang_CDC, 2015_Wang_H2} only work for a particular type of locality constraint (called the $d$-localized constraint), while the algorithm proposed here can handle arbitrary sparsity constraints.
\end{remark}

\section{Convergence and Stopping Criteria for Algorithm \eqref{eq:admm}} \label{sec:ADMM}
Assumption \ref{as:1} - \ref{as:3} imply feasibility and strong duality of optimization problem \eqref{eq:gopt-od} (or its equivalent formulation specified in equation \eqref{eq:sub3}). From Assumptions \ref{as:2} and \ref{as:3}, we know that the extended-real-value functionals $\rowdec{h}(\cdot)$ and $\coldec{h}(\cdot)$ defined in equation \eqref{eq:extended} are closed, proper, and convex. Under these assumptions, problem \eqref{eq:sub3} satisfies the convergence conditions stated in \cite{BoydADMM}. From \cite{BoydADMM}, we have objective convergence, dual variable convergence ($\ttf{\Lambda}^{k} \to \ttf{\Lambda}^{*}$ as $k \to \infty$), and residual convergence ($\ttf{\Phi}^k - \ttf{\Psi}^k \to 0$ as $k \to \infty$). 
Note that \eqref{eq:sub3} may not have a unique optimal point, so that the optimization variables $\ttf{\Phi}^k$ and $\ttf{\Psi}^k$ do not necessarily converge. If we further assume that the functional objective $g(\cdot)$ is strongly convex with respect to $\ttf{\Phi}$, then problem \eqref{eq:sub3} has a unique optimal solution $\ttf{\Phi}^*$. In this case, objective convergence implies primal variable convergence, and hence we have that $\ttf{\Phi}^k \to \ttf{\Phi}^*$ and $\ttf{\Psi}^k \to \ttf{\Phi}^*$ as $k \to \infty$.

The stopping criteria is designed using the methods suggested in \cite{BoydADMM}, wherein we use $|| \ttf{\Phi}^k - \ttf{\Psi}^k ||_{\mathcal{H}_2}$ as a measure of primal infeasibility and $|| \ttf{\Psi}^k - \ttf{\Psi}^{k-1} ||_{\mathcal{H}_2}$ as a measure of dual infeasibility. We further note that under the separability and locality conditions described in this paper, the square of these two functions can be computed in a localized and parallel way. The algorithm \eqref{eq:admm1} - \eqref{eq:admm3} terminates when $|| \ttf{\Phi}^k - \ttf{\Psi}^k ||_{\mathcal{H}_2} < \epsilon^{pri}$ and $|| \ttf{\Psi}^k - \ttf{\Psi}^{k-1} ||_{\mathcal{H}_2} < \epsilon^{dual}$ are satisfied for some feasibility tolerances $\epsilon^{pri}$ and $\epsilon^{dual}$.

In practice, we may not know \emph{a priori} whether optimization problem \eqref{eq:gopt-od} is feasible. In other words, we may not know whether Assumption \ref{as:1} holds. Consider the case that the ADMM subroutines \eqref{eq:admm1} - \eqref{eq:admm3} are solvable, but problem \eqref{eq:sub3} is infeasible. In this situation, the stopping criteria on primal infeasibility $|| \ttf{\Phi}^k - \ttf{\Psi}^k ||_{\mathcal{H}_2} < \epsilon^{pri}$ may not be satisfied. To avoid an infinite number of iterations, we set a limit on the number of iterations in the ADMM algorithm. In this way, whether the ADMM algorithm  as a scalable algorithm to check the feasibility of the SLS problem \eqref{eq:of_lop}.  
If the ADMM algorithm does not converge, then we know that SLS problem is not feasible with respect to the specified SLCs.
A general method for designing feasible SLCs is beyond the scope of this paper.  As an example of such an algorithm howeve, we refer the reader to our recent paper \cite{2015_Wang_Reg}, in which we present a method that allows for the joint design of an actuator architecture and a corresponding pair of feasible locality and FIR SLCs.

\section{Expressing ADMM Updates using Proximal Operators} \label{sec:appendix}
Here we explain how to express the solutions of to the iterate update subproblems \eqref{eq:admm1} and \eqref{eq:admm2} using proximal operators. We focus our discussion on subproblem \eqref{eq:admm2}, or its equivalent formulation described in equation \eqref{eq:admm_lop}, as an analogous argument holds for subproblem \eqref{eq:admm1}. Recall that the set $\coldec{\s}$ in problem \eqref{eq:admmlop3} is given by $\coldec{\s} = \Ell \cap \FT \cap \coldec{\Sother}$. We use the column-wise partition and the dimension reduction techniques described in Appendix \ref{sec:reduction} to simplify subproblem \eqref{eq:admm_lop} to
\begin{subequations} \label{eq:decom4}
\begin{align} 
\underset{\Sub{\ttf{\Psi}}{\rone_j}{\colps_j}}{\text{minimize }} \quad & \coldec{g}_j(\Sub{\ttf{\Psi}}{\rone_j}{\colps_j}) \nonumber\\
+ & \frac{\rho}{2} || \Sub{\ttf{\Psi}}{\rone_j}{\colps_j} - \Sub{\ttf{\Phi}^{k+1}}{\rone_j}{\colps_j} - \Sub{\ttf{\Lambda}^{k}}{\rone_j}{\colps_j} ||_{\mathcal{H}_2}^2 \label{eq:decom4-1}\\
\text{subject to } \quad & \Sub{\tf Z_{AB}}{\rtwo_j}{\rone_j} \Sub{\ttf{\Psi}}{\rone_j}{\colps_j} = \Sub{J_B}{\rtwo_j}{\colps_j} \label{eq:decom4-2} \\
& \Sub{\ttf{\Psi}}{\rone_j}{\colps_j} \in \Sub{\Ell}{\rone_j}{\colps_j} \cap \FT \cap \coldec{\Sother}_j \label{eq:decom4-3}
\end{align}
\end{subequations}
for $j = 1, \dots p$. 
In the simplified problem \eqref{eq:decom4}, the transfer matrices $\Sub{\ttf{\Psi}}{\rone_j}{\colps_j}$, $\Sub{\ttf{\Phi}^{k+1}}{\rone_j}{\colps_j}$, and $\Sub{\ttf{\Lambda}^{k}}{\rone_j}{\colps_j}$ all have a FIR with a horizon of $T$. We can therefore represent the optimization variables of problem \eqref{eq:decom4} as column vectors, i.e., as 
\begin{equation}
\Psi_{v(j)} = \text{vec}(\begin{bmatrix} \Sub{\Psi}{\rone_j}{\colps_j}[0] & \cdots & \Sub{\Psi}{\rone_j}{\colps_j}[T] \end{bmatrix}) \nonumber
\end{equation} 
where $\Psi_{v(j)}$ is the vectorization of all the spectral components of $\Sub{\ttf{\Psi}}{\rone_j}{\colps_j}$. Similarly, we define $\Phi_{v(j)}^{k+1}$ and $\Lambda_{v(j)}^k$ as the vectorization of all of the spectral components of $\Sub{\ttf{\Phi}^{k+1}}{\rone_j}{\colps_j}$ and $\Sub{\ttf{\Lambda}^{k}}{\rone_j}{\colps_j}$, respectively.
The optimization problem \eqref{eq:decom4} can then be written as
\begin{subequations} \label{eq:decom5}
\begin{align} 
\underset{\Psi_{v(j)}}{\text{minimize }} \quad & g_{v(j)}(\Psi_{v(j)}) + \frac{\rho}{2} || \Psi_{v(j)} - \Phi_{v(j)}^{k+1} - \Lambda_{v(j)}^k ||_2^2 \label{eq:decom5-1}\\
\text{subject to } \quad & \Psi_{v(j)} \in \s_{v(j)} \label{eq:decom5-2}
\end{align}
\end{subequations}
where $g_{v(j)}(\cdot)$ is the vectorized form of $\coldec{g}_j(\cdot)$, and $\s_{v(j)}$ is the set constraint imposed by \eqref{eq:decom4-2} - \eqref{eq:decom4-3}. Define the indicator function as
\begin{equation}
\mathcal{I}_{\mathcal{D}}(x) = \left\{ \begin{array}{ll}
0 & x \in \mathcal{D}\\
\infty & x \not \in \mathcal{D}.
\end{array} \right. \nonumber
\end{equation}
This allows us to rewrite \eqref{eq:decom5} as the following unconstrained problem 
\begin{align}
\underset{\Psi_{v(j)}}{\text{minimize }} \quad &\mathcal{I}_{\s_{v(j)}}(\Psi_{v(j)}) + g_{v(j)}(\Psi_{v(j)}) \nonumber\\
&+ \frac{\rho}{2} || \Psi_{v(j)} - \Phi_{v(j)}^{k+1} - \Lambda_{v(j)}^k ||_2^2. \label{eq:decom6}
\end{align}
The solution of the unconstrained problem \eqref{eq:decom6} can then be expressed using the proximal operator as
\begin{equation}
\Psi_{v(j)}^{k+1} = \textbf{prox}_{\mathcal{I}_{\s_{v(j)}} + \frac{1}{\rho} g_{v(j)}} ( \,\Phi_{v(j)}^{k+1} + \Lambda_{v(j)}^k \,). \label{eq:proximal}
\end{equation}
Equation \eqref{eq:proximal} is a solution of \eqref{eq:decom4}.
Therefore, the ADMM update \eqref{eq:admm2} can be carried out by evaluating the proximal operators \eqref{eq:proximal} for $j = 1, \dots p$, all in the reduced dimension.



\ifCLASSOPTIONcaptionsoff
  \newpage
\fi



\bibliographystyle{IEEEtran}
\bibliography{Distributed_new}
%
%
%

%

%
%
%

\begin{IEEEbiographynophoto}
{Yuh-Shyang Wang} is currently pursuing his PhD in Control and Dynamical Systems at Caltech. 
His research interests include optimal control, optimization, and system design.
\end{IEEEbiographynophoto}
%
\begin{IEEEbiographynophoto}
{Nikolai Matni} is a postdoctoral scholar in Computing and Mathematical Sciences at Caltech.
He received his PhD in Control and Dynamical Systems from Caltech.
His research interests broadly encompass the use of layering, dynamics, control and optimization in the design and analysis of complex cyber-physical systems; current application areas include software defined networking and sensorimotor control. 
His work on controller architecture design was awarded the 2013 IEEE CDC Best Student-Paper Award.
\end{IEEEbiographynophoto}
%
\begin{IEEEbiographynophoto}{John C. Doyle} is Professor of Control and Dynamical Systems, Electrical Engineer, and BioEngineering at Caltech. BS, MS EE, MIT (1977), PhD, Math, UC Berkeley (1984). Research is on mathematical foundations for complex networks with applications in biology, technology, medicine, ecology, and neuroscience. Paper prizes include IEEE Baker, IEEE Automatic Control Transactions (twice), ACM Sigcomm, and ACC American Control Conference. Individual awards include AACC Eckman, IEEE Control Systems Field, and IEEE Centennial Outstanding Young Engineer Awards. Has held national and world records and championships in various sports.
\end{IEEEbiographynophoto}




\end{document}